# FINITELY ADDITIVE BELIEFS AND UNIVERSAL TYPE SPACES[1]

By Martin Meier

*Instituto de Análisis Económico, CSIC*

The probabilistic type spaces in the sense of Harsanyi [*Management Sci.* **14** (1967/68) 159–182, 320–334, 486–502] are the prevalent models used to describe interactive uncertainty. In this paper we examine the existence of a universal type space when beliefs are described by finitely additive probability measures. We find that in the category of all type spaces that satisfy certain measurability conditions ($\kappa$-measurability, for some fixed regular cardinal $\kappa$), there is a universal type space (i.e., a terminal object) to which every type space can be mapped in a unique beliefs-preserving way. However, by a probabilistic adaption of the elegant sober-drunk example of Heifetz and Samet [*Games Econom. Behav.* **22** (1998) 260–273] we show that if all subsets of the spaces are required to be measurable, then there is no universal type space.

**1. Introduction.** Consider players that are uncertain about a set $S$, called the *space of states of nature,* each element of which can be thought of as a complete list of the players' strategy sets and payoff functions, that is, a complete specification of the "rules" of the game that depend on the state of nature. (Other interpretations are also possible. For example, if a game of complete information is given, a state $s \in S$ could be the strategy profile that the players are actually going to choose; see the analysis of epistemic conditions for Nash equilibrium by Aumann and Brandenburger [1].) In such a situation, following a Bayesian approach, each player will base his choice of a strategy on his subjective beliefs (i.e., a probability measure) on $S$.

---

Received May 2002; revised May 2005.

[1]Supported by the DFG (German Science Foundation) via the Graduiertenkolleg "Mathematische Wirtschaftsforschung" (University of Bielefeld), the European Union via the TMR-Network "Cooperation and Information" (University of Tel Aviv and Université de Caen) and a Marie Curie Individual Fellowship at CORE (Université Catholique de Louvain), and the Spanish Ministerio de Educación y Ciencia via a Ramon y Cajal Fellowship (IAE-CSIC) and Research Grant SEJ2004-07861.

Since a player's payoff depends also on the choices of the other players, and these are based on their beliefs as well, each player must also have beliefs on the other players' beliefs on $S$. For the same reason, he must also have beliefs on other players' beliefs on his beliefs on $S$, beliefs on other players' beliefs on his beliefs on their beliefs on $S$, and so on. So, in analyzing such a situation, it seems to be unavoidable to work with infinite hierarchies of beliefs. Thus, the resulting model is complicated and cumbersome to handle. In fact, this was the issue that prevented for a long time the analysis of games of incomplete information.

A major breakthrough took place with three articles of Harsanyi [5], where he succeeded in finding another, more workable model to describe interactive uncertainty. He invented the notions of *type* and *type space*: With each point in a type space, called a *state of the world*, are associated a state of nature and, for each player, a probability measure on the type space itself (i.e., that player's *type* in this state of the world). Usually it is assumed that the players "know their own type," that is, a type of a player in a state assigns probability 1 to the set of those states where this player is of this type. This is the formalization of the idea that the players should be self-conscious. Since each state of the world is associated with a state of nature, each player's type in a state of the world induces a probability measure on $S$. But also, since with each state of the world there is associated a type for each player (and hence indirectly a probability measure on $S$ for this player), the type of a player in a state of the world induces a probability measure on the other players' probability measures on $S$. Proceeding like this, one obtains in each state of the world a hierarchy of beliefs for each player, in the sense described above.

The advantages of Harsanyi's model are obvious: Since we have in each state of the world just one probability measure for each player, contrary to the hierarchical description of beliefs, this model fits in the classical Bayesian framework of describing beliefs by one probability measure, and provides therefore all its advantages (e.g., it allows for integration with respect to beliefs).

However, there are also several serious questions that arise with the use of this model: Although each state of the world in a type space induces a hierarchy of beliefs for each player, the converse is not obvious: Does each profile of hierarchies of beliefs arise from a state of the world in some type space, and if so, is there a type space such that every profile of belief hierarchies is generated by some state in this type space? What "are" the states of the world, and what justifies using one particular type space and not another? In particular, contrary to the case of the hierarchical description of beliefs, it is not clear what "all possible types" (resp. "all possible states of the world") are. More precisely, given a game of incomplete information, working in one fixed type space to analyze that game could be restrictive in the sense that we might miss some possible types (resp. possible states of the world) that are just present in a bigger type space that contains the one we use. If this were the case for every type space, then the use of type spaces would be problematic from a theoretical point of view, because of the restrictive character of this concept, but it would be problematic also from a more practical point of view: In their contributions to the debate on epistemic conditions for backward induction in perfect information games, Stalnaker [16] and Battigalli and Siniscalchi [2] have pointed out that the players do "their best" to rationalize their opponents' behavior if the backward-induction outcome is to obtain. This translates into using a type space where a player can find the needed types he has to attribute to the other players, if he has to explain (i.e., rationalize) the others' behavior.

The question concerning "all possible types" can be answered and the related problems can be solved if there is a type space to which every type space (on the same space of states of nature and for the same set of players, of course) can be mapped, preferably always in a unique way, by a map that preserves the structure of the type space, that is, the manner in which types and states of nature are associated with states of the world, a so-called *type morphisms*. Such a type space would be called a *universal type space*. If such a space always exists, one could, in principle, carry out the analysis of a game of incomplete information in the corresponding universal type space without any risk of missing a relevant state of the world. On a technical level, the type spaces—on a fixed set of states of nature and for a fixed player set—as objects and the type morphisms as morphisms form a category. If we always require the map from a type space to the universal type space to be unique, then, if it exists, such a universal type space is a *terminal object* of this category. A terminal object of a category is known to be unique up to isomorphism. Hence, we are justified to talk about *the* universal type space.

The existence of a universal type space was proved by Mertens and Zamir [14] under the assumption that the underlying space of states of nature is a compact Hausdorff space and all involved functions are continuous. That topological assumption was relaxed by Brandenburger and Dekel [3], Heifetz [6] and Mertens, Sorin and Zamir [13] to more general topological assumptions. Finally, the general measure-theoretic case was solved by Heifetz and Samet [10], who showed that there also exists a universal type space in this case. However, in all these articles it has always been assumed that the players' beliefs are $\sigma$-additive. This seems to be a rather strong assumption on the epistemic attitudes of the players.

Savage's postulates [15] imply subjective probabilities that are finitely but not countably additive. Given the importance of Savage's theory within decision theory (i.e., "one-player game theory"), it is natural to ask—and desirable to know—what happens if we describe the beliefs of players accordingly in an interactive context (games) by finitely additive probability measures? Does there still exist a universal type space?

Still, we are dealing with (now finitely additive) measures, so we have to define a field of events. The question arises as to which measurability condition is the right one. Should the field of events be just a field, a $\sigma$-field or should we assume that all subsets of the space of states of the world are events, that is, that the field is simply the power set? As this question does not seem to have a clear-cut answer, we analyze the existence/nonexistence of a universal type space for finitely additive beliefs for a whole class of different measurability conditions that include the three above-mentioned cases.

We introduce $\kappa$-fields, where $\kappa$ is an (usually regular) infinite cardinal number, as fields that are closed under the intersection of every set of events (i.e., subset of the field) that has cardinality strictly less than $\kappa$. It follows that $\aleph_0$-fields are fields in the usual sense and $\aleph_1$-fields are $\sigma$-fields in the usual sense. Then, we define $\infty$-fields as fields that are closed under the intersection of every set of events (of cardinality whatsoever).

We define $\kappa$-type spaces as type spaces where the set of measurable events in the set of states of the world, as well as the set of measurable events in the set of states of nature, is a $\kappa$-field, and $\infty$-type spaces as type spaces where the set of measurable events in the set of states of nature is the full power set and the set of measurable events in the set of states of the world is a $\infty$-field. Also, we define $*$-type spaces as type spaces where the set of measurable events in the set of states of the world as well as the set of measurable events in the set of states of nature is the full power set. Furthermore, we define type morphisms, that is, structure-preserving maps from one $\kappa$-type space (resp. $\infty$-type space or $*$-type space) to another (not necessarily different) one.

Given a nonempty set of players $I$, a nonempty set of states of nature $S$ and a $\kappa$-field $\Sigma_S$ on $S$, we define, similar to Heifetz and Samet [10], a kind of modal language, the formulas of which we call $\kappa$-expressions. But if $\kappa$ is uncountable, contrary to Heifetz and Samet [10], we allow also for formulas of infinite length (but strictly less than $\kappa$). Then, we collect the $\kappa$-descriptions (by means of $\kappa$-expressions) of all states of the world in all $\kappa$-type spaces on $S$ for player set $I$. Then we show that the set of $\kappa$-descriptions can be endowed with the structure of a $\kappa$-type space (Proposition 4). In this way, we construct in Section 4 a universal $\kappa$-type space on $S$ for player set $I$ to which every $\kappa$-type space on $S$ for player set $I$ can be mapped by a unique type morphism (Theorem 1).

As Heifetz [7] has shown, there are consistent hierarchies of finitely additive (in fact even $\sigma$-additive) beliefs up to—but excluding—level $\omega$ (i.e., the first infinite level), that have at least two different finitely additive extensions to level $\omega$. Does a similar phenomenon hold also on the higher transfinite levels of consistent hierarchies? Or, put differently in terms of expressions rather than hierarchies, is there, on the contrary, a regular cardinal $\widehat{\kappa}$ such that for all (regular) cardinals $\kappa > \widehat{\kappa}$, the $\widehat{\kappa}$-description of a state in a $\kappa$-type space determines already the $\kappa$-description? If this were the case, it would be unnecessary to consider $\kappa$-type spaces for $\kappa > \widehat{\kappa}$ and we could restrict ourselves to $\kappa$-type spaces for $\kappa \leq \widehat{\kappa}$. We show in Theorem 3, by using a probabilistic adaptation of the "sober-drunk" example of Heifetz and Samet [9], that—with at least two players and two states of nature—this is not the case. Hence, it makes sense to consider $\kappa$-type spaces for every (regular) infinite cardinal $\kappa$.

Also, this example implies that—again, with at least two players and two states of nature—there is no universal $\infty$-type space and no universal $*$-type space (Theorem 4 and Corollary 1), even if we do not require the morphisms from the type spaces to the universal type space to be unique.

**2. Preliminaries.** First, we will define $\kappa$-measurable spaces, where $\kappa$ denotes here an (usually regular) infinite cardinal number. For an introduction to ordinal and cardinal numbers, see [4] or any other textbook on set theory. Then, we will develop parts of the theory of $\kappa$-measurable spaces needed in the sequel, collect some known facts about finitely additive (probability) measures and define the main objects of our study in this paper, the $\kappa$-, $\infty$- and $*$-type spaces.

An infinite cardinal number $\kappa$ is called *regular*, if it is not the supremum of a set of less than $\kappa$-many ordinal numbers which are all strictly smaller than $\kappa$. For example, $\aleph_0$ and all the $\aleph_{\alpha+1}$ are regular, while $\aleph_\omega$ is *singular* (i.e., infinite and not regular) ($\aleph_\omega = \sup\{\aleph_n | n < \omega\}$), where $\omega$ denotes here the first infinite ordinal number. For a set $M$, denote by $|M|$ the cardinality of $M$.

Unless otherwise stated, $\alpha, \beta, \gamma, \zeta, \eta, \xi$ denote ordinal numbers, $\delta$ delta-measures, $\theta$ functions from the set of states of the world to the set of states of nature, $\kappa$ cardinal numbers, $\lambda$ limit ordinal numbers, $\mu$ and $\nu$ measures, $\pi$ projections, $\varphi, \chi, \psi$ expressions and $\omega$, apart from above, sets of expressions. For a set $M$, $\mathrm{Pow}(M)$ denotes the set of all subsets of $M$, that is, the power set.

Let $\kappa$ be an infinite cardinal number and $M$ a nonempty set. A $\kappa$-*field on* $M$ is a field $\Sigma$ on $M$ such that $\mathcal{E} \subseteq \Sigma$ and $|\mathcal{E}| < \kappa$ imply $\bigcap \mathcal{E} := \bigcap_{E \in \mathcal{E}} E \in \Sigma$. It follows that $\mathcal{E} \subseteq \Sigma$ and $|\mathcal{E}| < \kappa$ imply $\bigcup \mathcal{E} := \bigcup_{E \in \mathcal{E}} E \in \Sigma$.

Consequently, a $\kappa$-*measurable space* is a pair $(M, \Sigma)$, where $M$ is a nonempty set and $\Sigma$ is a $\kappa$-field on $M$.

A set of subsets of a nonempty set is a $\aleph_0$-field iff it is a field and it is a $\aleph_1$-field iff it is a $\sigma$-field. If $\kappa' < \kappa$, then every $\kappa$-field is also a $\kappa'$-field.

REMARK 1. *Let $\kappa$ be a singular cardinal number and $(M, \Sigma)$ be a $\kappa$-measurable space. Then $\Sigma$ is already a $\kappa^+$-field, where $\kappa^+$ denotes the successor cardinal of $\kappa$.*

PROOF. Let $\mathcal{E} \subseteq \Sigma$ such that $|\mathcal{E}| \leq \kappa$. So, $\mathcal{E}$ has the form $\{E_\alpha | \alpha < \kappa\}$. Let $\widehat{\kappa} < \kappa$ be the cofinality of $\kappa$. Then there is a function $f : \widehat{\kappa} \to \kappa$, such that $\bigcup_{\beta < \widehat{\kappa}} f(\beta) = \kappa$. Note that $|f(\beta)| < \kappa$, for $\beta < \widehat{\kappa}$. It follows that $\bigcap_{\alpha < \kappa} E_\alpha = \bigcap_{\beta < \widehat{\kappa}} (\bigcap_{\alpha < f(\beta)} E_\alpha) \in \Sigma$. Since $\Sigma$ is a field, it follows that it is a $\kappa^+$-field. $\square$

Since $\kappa^+$ is always regular, the above remark shows that it is redundant to consider $\kappa$-fields ($\kappa$-measurable spaces, resp.) if $\kappa$ is a singular cardinal.

Let $M$ be a nonempty set. A $\infty$-*field on $M$* is a field $\Sigma$ on $M$ such that $\mathcal{E} \subseteq \Sigma$ implies $\bigcap \mathcal{E} := \bigcap_{E \in \mathcal{E}} E \in \Sigma$. Again, it follows that $\mathcal{E} \subseteq \Sigma$ implies $\bigcup \mathcal{E} := \bigcup_{E \in \mathcal{E}} E \in \Sigma$.

Accordingly, a $\infty$-*measurable space* is a pair $(M, \Sigma)$, where $M$ is a nonempty set and $\Sigma$ is a $\infty$-field on $M$.

A $*$-*measurable space* is a pair $(M, \mathrm{Pow}(M))$, where $M$ is a nonempty set.

Note that every $*$-measurable space is a $\infty$-measurable space and every $\infty$-measurable space is a $\kappa$-measurable space for every infinite ordinal $\kappa$. A $\infty$-measurable space $(M, \Sigma)$ is a $*$-measurable space iff for all $m \neq m' \in M$ there is an $E \in \Sigma$ such that $m \in E$ and $m' \notin E$.

EXAMPLE 1. Let $M = \{0, 1\}$, $\Sigma = \{\varnothing, M\}$. $(M, \Sigma)$ is $\infty$-measurable, but not $*$-measurable.

DEFINITION 1. Let $M$ be a nonempty set and $\mathcal{F}$ a field on $M$. A *finitely additive measure on* $(M, \mathcal{F})$ is a function $\mu : \mathcal{F} \to \mathbb{R} \cup \{+\infty\}$, such that:

(i) $0 \leq \mu(F)$, for all $F \in \mathcal{F}$,

(ii) $\mu(E \cup F) = \mu(E) + \mu(F)$, for all disjoint $E, F \in \mathcal{F}$.

$\mu$ is a *finitely additive probability measure on* $(M, \mathcal{F})$, if in addition

(iii) $\mu(M) = 1$.

DEFINITION 2. Let $M$ be a nonempty set, $\mathcal{F}$ a field on $M$, $\mu$ a finitely additive measure on $(M, \mathcal{F})$, and $E \subseteq M$.

We define the *outer measure of $E$ induced by* $\mu$ as

$$\mu^*(E) := \inf\{\mu(F) | F \in \mathcal{F} \text{ such that } E \subseteq F\},$$

and the *inner measure of $E$ induced by* $\mu$ as

$$\mu_*(E) := \sup\{\mu(F) | F \in \mathcal{F} \text{ such that } F \subseteq E\}.$$

If not stated otherwise, we keep the following.

CONVENTION 1. (i) If $(M, \Sigma)$ is a $\kappa$-measurable space, then $\Delta^\kappa(M, \Sigma)$ denotes the space of finitely additive probability measures on $(M, \Sigma)$. We consider this space itself as a $\kappa$-measurable space endowed with the $\kappa$-field $\Sigma_{\Delta^\kappa}$ generated by all the sets $\{\mu \in \Delta^\kappa(M, \Sigma) | \mu(E) \geq p\}$, where $E \in \Sigma$ and $p \in [0, 1]$.

(ii) Similarly, we denote by $\Delta(M, \mathrm{Pow}(M))$ the set of all finitely additive probability measures on $(M, \mathrm{Pow}(M))$.

Of course, (i) of this convention depends on the particular $\kappa$ chosen.

REMARK 2. *Let $(M', \Sigma')$ and $(M, \Sigma)$ be $\kappa$-measurable spaces and let $f: M' \to M$ be measurable. Then:*

(a) *If $\mu'$ is a finitely additive probability measure on $(M', \Sigma')$, then $\mu'(f^{-1}(\cdot))$ [i.e., $\mu'(f^{-1}(E))$, for $E \in \Sigma$] is a finitely additive probability measure on $(M, \Sigma)$.*

(b) *If $\Delta_f^\kappa: \Delta^\kappa(M', \Sigma') \to \Delta^\kappa(M, \Sigma)$ is defined by $\Delta_f^\kappa(\mu') := \mu'(f^{-1}(\cdot))$, for $\mu' \in \Delta^\kappa(M', \Sigma')$, then $\Delta_f^\kappa$ is measurable, since we have $\Delta_f^\kappa(\mu')(E) \geq p$ iff $\mu'(f^{-1}(E)) \geq p$, for $E \in \Sigma$.*

REMARK 3. *Let $(M', \Sigma')$ and $(M, \Sigma)$ be $\kappa$-measurable spaces and let $f: M' \to M$ be measurable and onto. Then:*

(a) *$f^{-1}(\Sigma) := \{f^{-1}(E) | E \in \Sigma\}$ is a $\kappa$-field on $M'$ and a subset of $\Sigma'$.*

(b) *If $\mu$ is a finitely additive measure on $(M, \Sigma)$, then $\mu$ induces a finitely additive measure $\mu'$ on $(M', f^{-1}(\Sigma))$ defined by $\mu'(f^{-1}(E)) := \mu(E)$. Furthermore, if $\mu$ is a finitely additive probability measure, then $\mu'$ is a finitely additive probability measure.*

LEMMA 1. *Let $\gamma < \alpha$ be ordinal numbers. For $\gamma \leq \beta < \alpha$ let $(M^\beta, \mathcal{F}^\beta)$ be a $\aleph_0$-measurable space (i.e., $M^\beta$ is a nonempty set and $\mathcal{F}^\beta$ is a field on $M^\beta$) and $\mu^\beta$ a finitely additive probability measure on $(M^\beta, \mathcal{F}^\beta)$, let $M^\alpha$ be a nonempty set, and for $\gamma \leq \xi < \zeta \leq \alpha$ let $f_{\xi, \zeta}: M^\zeta \to M^\xi$ be onto and, if $\zeta < \alpha$, let $f_{\xi, \zeta}$ be $\mathcal{F}^\zeta - \mathcal{F}^\xi$-measurable, such that:*

1. *$f_{\xi, \beta} \circ f_{\beta, \zeta} = f_{\xi, \zeta}$, for all $\xi < \beta < \zeta$ such that $\gamma \leq \xi < \beta < \zeta \leq \alpha$,*
2. *$\mu^\beta(f_{\xi, \beta}^{-1}(E^\xi)) = \mu^\xi(E^\xi)$, for all $\xi < \beta$ such that $\gamma \leq \xi < \beta < \alpha$ and all $E^\xi \in \mathcal{F}^\xi$.*

*Then:*

(a) *$\bigcup_{\gamma \leq \beta < \alpha} f_{\beta, \alpha}^{-1}(\mathcal{F}^\beta)$ is a field on $M^\alpha$,*

(b) *$(\mu^\beta)_{\gamma \leq \beta < \alpha}$ induces a well-defined finitely additive probability measure $\mu^{<\alpha}$ on $(M^\alpha, \bigcup_{\gamma \leq \beta < \alpha} f_{\beta, \alpha}^{-1}(\mathcal{F}^\beta))$, defined by $\mu^{<\alpha}(f_{\beta, \alpha}^{-1}(E^\beta)) := \mu^\beta(E^\beta)$, for $E^\beta \in \mathcal{F}^\beta$.*

PROOF. That $\mu^{<\alpha}$ is well-defined follows from the above conditions 1 and 2 and the fact that the $f_{\beta,\alpha}$'s are onto. In light of the preceding remark, the rest is clear. □

NOTATION 1. Let $M$ be a nonempty set, $\mathcal{F}$ a field on $M$, and $E \subseteq M$. Then denote by $[\mathcal{F}, E]$ the set of all subsets of $M$ of the form $(L \cap E) \cup (N \cap (M \setminus E))$, where $L, N \in \mathcal{F}$. It is easy to check that $[\mathcal{F}, E]$ is the smallest field that extends $\mathcal{F}$ and contains $E$ as an element.

For further reference, we cite the following two lemmas (in a somewhat different form), which are theorems by Łoś and Marczewski [12] and Horn and Tarski [11].

LEMMA 2. Let $M$ be a nonempty set, $\mathcal{F}$ a field on $M$, $E \subseteq M$, $\mu$ a finitely additive probability measure on $(M, \mathcal{F})$, $\mu_*(E)$ the inner measure of $E$, $\mu^*(E)$ the outer measure of $E$, and $p$ a real number such that $\mu_*(E) \leq p \leq \mu^*(E)$.
  Then there exists a finitely additive probability measure $\nu$ that extends $\mu$ to the field $[\mathcal{F}, E]$ such that $\nu(E) = p$.

PROOF. Follows directly from Theorem 2 of [12]. □

Sometimes, we will refer to the above lemma as the "Łoś–Marczewski theorem."

LEMMA 3. Let $\mathcal{F}_1 \subseteq \mathcal{F}_2$ be fields on the nonempty set $M$ and let $\mu$ be a finitely additive probability measure on $(M, \mathcal{F}_1)$. Then there exists an extension of $\mu$ to a finitely additive probability measure $\nu$ on $(M, \mathcal{F}_2)$.

PROOF. Follows from point (i) of Section 4 of [12] and also from [11], page 477, Theorem 1.21. □

**3. Type spaces.** For this section, unless otherwise stated, we fix a regular cardinal $\kappa$. Furthermore we fix a nonempty set of players $I$, a nonempty set of states of nature $S$, and, unless otherwise stated, a $\kappa$-field $\Sigma_S$ on $S$, such that for all $s, s' \in S$ with $s \neq s'$ there is an $E \in \Sigma_S$ such that $s \in E$ and $s' \notin E$.

We define now $\kappa$-type spaces, $\infty$-type spaces and $*$-type spaces, that is, the objects which we will study in this paper.

DEFINITION 3. A $\kappa$-type space on $S$ for player set $I$ is a 4-tuple

$$\underline{M} := \langle M, \Sigma, (T_i)_{i \in I}, \theta \rangle,$$

where:

(a) $M$ is a nonempty set,

(b) $\Sigma$ is a $\kappa$-field on $M$,

(c) for $i \in I$: $T_i$ is a $\Sigma - \Sigma_{\Delta^\kappa}$-measurable function from $M$ to $\Delta^\kappa(M, \Sigma)$, the space of finitely additive probability measures on $(M, \Sigma)$, such that for all $m \in M$ and $A \in \Sigma$: $[T_i(m)] \subseteq A$ implies $T_i(m)(A) = 1$, where $[T_i(m)] := \{m' \in M | T_i(m') = T_i(m)\}$,

(d) $\theta$ is a $\Sigma - \Sigma_S$-measurable function from $M$ to $S$.

This structure is interpreted as follows: $M$ is the set of states of the world. Such a state determines completely the objective parameters of the players' interaction, that is, the state of the nature $\theta(m)$, as well as the players' beliefs about the true state of the world. In general, in a state of the world $m \in M$, player $i$ will not know the true state of the world $m$; he will just have a probability measure $T_i(m)$ over the set of states of the world. $T_i(m)$ describes his beliefs in state $m$, that is, the *type of player $i$ in state $m$*. [Knowing $m$ would mean that $T_i(m) = \delta_m$, where $\delta$ denotes the Kronecker delta.] $\theta(m)$ is the state of nature that corresponds, to the state of the world $m$. While there might be many states of the world to which a given state of nature $s \in S$ corresponds, we have that to every state of the world there corresponds one and only one state of nature. This is expressed by the fact that $\theta$ is a function $\theta : M \to S$.

We will refer to the property that for all $i \in I$, $m \in M$ and $A \in \Sigma$: $[T_i(m)] \subseteq A$ implies $T_i(m)(A) = 1$, as the *introspection property* of $\kappa$-type spaces ($\infty$-type spaces and $*$-type spaces, resp., see below). This expresses the self-consciousness of the players: In a state of the world $m$ a player does not attribute a positive probability to states where he has a different belief from the belief he has in the present state $m$.

Doing obvious changes, the proofs (of the theorems in this paper) would go through, if we were to abandon this property; in fact, things would be easier then.

DEFINITION 4.  A *$\infty$-type space on $S$ for player set $I$ is a 4-tuple*

$$\underline{M} := \langle M, \Sigma, (T_i)_{i \in I}, \theta \rangle,$$

where:

(a) $M$ is a nonempty set,

(b) $\Sigma$ is a $\infty$-field on $M$,

(c) for $i \in I$: $T_i$ is a measurable function from $(M, \Sigma)$ to $\Delta^\infty(M, \Sigma)$, the space of finitely additive probability measures on $(M, \Sigma)$, endowed with the $\infty$-field generated by all the sets $\{\mu \in \Delta^\infty(M, \Sigma) | \mu(E) \geq p\}$, where $E \in \Sigma$ and $p \in [0, 1]$, such that for all $m \in M$ and $A \in \Sigma$: $[T_i(m)] \subseteq A$ implies $T_i(m)(A) = 1$, where $[T_i(m)] := \{m' \in M | T_i(m') = T_i(m)\}$,

(d) $\theta$ is a $\Sigma - \mathrm{Pow}(S)$-measurable function from $M$ to $S$.

Note that for $\mu \neq \nu \in \Delta^\infty(M, \Sigma)$ there is an $E \in \Sigma$ and a $p \in [0,1]$ such that $\mu(E) \geq p$ and $\nu(E) < p$. This implies that the $\infty$-field of $\Delta^\infty(M, \Sigma)$ is in fact $\mathrm{Pow}(\Delta^\infty(M, \Sigma))$, the full power set. Hence, by the measurability of $T_i$, we have $[T_i(m)] \in \Sigma$. So, in fact, the condition that $[T_i(m)] \subseteq A$ implies $T_i(m)(A) = 1$ reduces to $T_i(m)([T_i(m)]) = 1$.

By the definitions, it is obvious that every $\infty$-type space on $S$ is a $\kappa$-type space on $S$, for every regular $\kappa$. [Set $\Sigma_S := \mathrm{Pow}(S)$ in the $\kappa$-type space.]

DEFINITION 5.   A *-type space on $S$ for player set $I$ is a 3-tuple

$$\underline{M} := \langle M, (T_i)_{i \in I}, \theta \rangle,$$

where:

(a) $M$ is a nonempty set,

(b) for $i \in I : T_i$ is a function from $M$ to $\Delta(M, \mathrm{Pow}(M))$, the space of finitely additive probability measures on $(M, \mathrm{Pow}(M))$, such that for all $m \in M : T_i(m)([T_i(m)]) = 1$, where $[T_i(m)] := \{m' \in M | T_i(m') = T_i(m)\}$,

(c) $\theta$ is a function from $M$ to $S$.

Equivalently, a *-type space $\underline{M}$ can be written as $\underline{M} = \langle M, \Sigma, (T_i)_{i \in I}, \theta \rangle$, where $\Sigma = \mathrm{Pow}(M)$. So, every *-type space on $S$ is a $\infty$-type space on $S$. And therefore, it is also a $\kappa$-type space on $S$.

We define now the beliefs-preserving maps between type spaces.

DEFINITION 6.   Let $\underline{M}' = \langle M', \Sigma', (T_i')_{i \in I}, \theta' \rangle$ and $\underline{M} = \langle M, \Sigma, (T_i)_{i \in I}, \theta \rangle$ be $\kappa$-type spaces ($\infty$-type spaces, *-type spaces, resp.) on $S$ for player set $I$.

A function $f : M' \to M$ is a type morphism if it satisfies the following conditions:

1. $f$ is $\Sigma' - \Sigma$-measurable,
2. for all $m' \in M'$:

$$\theta'(m') = \theta(f(m')),$$

3. for all $m' \in M'$, $E \in \Sigma$ and $i \in I$:

$$T_i(f(m'))(E) = T_i'(m')(f^{-1}(E)).$$

Note that the above definition of a type morphism does not depend on $\kappa$; that is, if $\kappa < \kappa'$, and $\underline{M}$ and $\underline{M}'$ are $\kappa'$-type spaces ($\infty$-type spaces, *-type spaces, resp.), then $f : M' \to M$ is a type morphism from $\underline{M}'$ to $\underline{M}$ viewed as $\kappa'$-type spaces ($\infty$-type spaces, *-type spaces, resp.) iff it is a type morphism from $\underline{M}'$ to $\underline{M}$ viewed as $\kappa$-type spaces. Similarly, if $\underline{M}'$ and $\underline{M}$ are $*$-type spaces, then $f : M' \to M$ is a type morphism from $\underline{M}'$ to $\underline{M}$ viewed as $*$-type spaces iff it is a type morphism from $\underline{M}'$ to $\underline{M}$ viewed as $\infty$-type spaces. (Note that in the case of $*$-type spaces, every function $f : M' \to M$ is measurable.)

DEFINITION 7. A type morphism is a *type isomorphism*, if it is one-to-one, onto, and the inverse function is also a type morphism.

It is easy to see that a function $f : M' \to M$ is a type isomorphism iff it is a type morphism and isomorphism of the measurable spaces $(M', \Sigma')$ and $(M, \Sigma)$.

An easy check shows:

REMARK 4. *$\kappa$-type spaces on $S$ for player set $I$ ($\infty$-type spaces, $*$-type spaces, resp.), as objects, and type morphisms, as morphisms, form a category.*

DEFINITION 8. (i) A $\kappa$-type space $\underline{\Omega}$ on $S$ for player set $I$ ($\infty$-type space, $*$-type space, resp.) is *weak-universal* if for every $\kappa$-type space $\underline{M}$ on $S$ for player set $I$ ($\infty$-type space, $*$-type space, resp.) there is a type morphism from $\underline{M}$ to $\underline{\Omega}$.

(ii) A $\kappa$-type space $\underline{\Omega}$ on $S$ for player set $I$ ($\infty$-type space, $*$-type space, resp.) is *universal* if for every $\kappa$-type space $\underline{M}$ on $S$ for player set $I$ ($\infty$-type space, $*$-type space, resp.) there is a *unique* type morphism from $\underline{M}$ to $\underline{\Omega}$.

To keep the terms of the already existing type space literature, we use the term "universal type space," although in the language of category theory the term "terminal type space" would be more adequate, since the universal type space is a terminal object in the category of type spaces.

PROPOSITION 1. *If they exist, universal $\kappa$-type spaces on $S$ for player set $I$ ($\infty$-type spaces, $*$-type spaces, resp.) are unique up to type isomorphism.*

PROOF. If $\underline{\Omega}$ and $\underline{U}$ are universal $\kappa$-type spaces ($\infty$-type spaces, $*$-type spaces, resp.) (on the same space of states of nature and for the same player set, of course), then there are type morphisms $f : U \to \Omega$ and $g : \Omega \to U$. It is easy to check that the composite of two type morphisms is also a type morphism and that the identity is always a type morphism from a $\kappa$-type space $\underline{\Omega}$ ($\infty$-type space, $*$-type space, resp.) to itself. By the uniqueness, it follows that $g \circ f = \mathrm{id}_U$ and therefore $f$ is one-to-one and $g$ is onto, and $f \circ g = \mathrm{id}_\Omega$ and therefore $g$ is one-to-one and $f$ is onto. $f$ and $g$ are type morphisms and $f = g^{-1}$ and $g = f^{-1}$. $\square$

To prove the existence of a universal $\kappa$-type space on $S$ for player set $I$ is the goal of the next section.

**4. The universal $\kappa$-type space in terms of expressions.** Again, for this section, unless otherwise stated, we fix a regular cardinal $\kappa$, a nonempty player set $I$ and a $\kappa$-measurable space of states of nature $(S, \Sigma_S)$ such that for all $s, s' \in S$ with $s \neq s'$ there is an $E \in \Sigma_S$ such that $s \in E$ and $s' \notin E$.

Given these data, we define $\kappa$-*expressions* (allowing also for infinite conjunctions) which are natural generalizations of the *expressions* defined by Heifetz and Samet [10]. These are formulas that describe events defined solely in terms of nature and the players' beliefs. Expressions are defined in a similar fashion as, for example, the formulas of the probability logic of Heifetz and Mongin [8]. Analogous to Heifetz and Samet [10], given a $\kappa$-type space on $S$ for player set $I$ and a state of the world in this type space, we define the $\kappa$-*description* of this state as the set of those $\kappa$-*expressions* that are true in this state of the world. Then, we show that the set of all $\kappa$-descriptions constitutes a $\kappa$-type space (Proposition 4) and that this $\kappa$-type space is the universal $\kappa$-type space (Theorem 1).

DEFINITION 9. For a $\kappa$-type space $\langle M, \Sigma, (T_i)_{i \in I}, \theta \rangle$ on $S$ for player set $I$, $i \in I$, $E \in \Sigma$ and $p \in [0, 1]$ define

$$\overline{B_i^p}(E) := \{ m \in M \mid T_i(m)(E) \geq p \}.$$

Note that $\overline{B_i^p}(E) = T_i^{-1}(\{ \mu \in \Delta^\kappa(M, \Sigma) \mid \mu(E) \geq p \})$ and that $\overline{B_i^p}(E) \in \Sigma$, if $E \in \Sigma$.

DEFINITION 10. Given a $\kappa$-measurable space of states of nature $(S, \Sigma_S)$ and a nonempty player set $I$, the set $\Phi^\kappa$ of $\kappa$-*expressions* is the least set such that:

1. every $E \in \Sigma_S$ is a $\kappa$-expression,
2. if $\varphi$ is a $\kappa$-expression, then $\neg\varphi$ is a $\kappa$-expression,
3. if $\varphi$ is a $\kappa$-expression, then $B_i^p(\varphi)$ is a $\kappa$-expression, for $i \in I$ and $p \in [0, 1]$,
4. if $\Psi$ is a nonempty set of $\kappa$-expressions with $|\Psi| < \kappa$, then $\bigwedge_{\varphi \in \Psi} \varphi$ is a $\kappa$-expression.

If $\Psi$ is a nonempty set of $\kappa$-expressions with $|\Psi| < \kappa$, then we set $\bigvee_{\varphi \in \Psi} \varphi := \neg \bigwedge_{\varphi \in \Psi} \neg\varphi$.

Since we work here with a fixed regular $\kappa$, we omit sometimes the superscript $\kappa$.

DEFINITION 11. Let $\underline{M} := \langle M, \Sigma, (T_i)_{i \in I}, \theta \rangle$ be a $\kappa$-type space on $S$ for player set $I$. Define:

1. $E^{\underline{M}} := \theta^{-1}(E)$, for $E \in \Sigma_S$,
2. $(\neg\varphi)^{\underline{M}} := M \setminus \varphi^{\underline{M}}$, for $\varphi \in \Phi^\kappa$,

3. $(B_i^p(\varphi))^{\underline{M}} := \overline{B_i^p}(\varphi^{\underline{M}})$, for $\varphi \in \Phi^\kappa$, $i \in I$ and $p \in [0,1]$,

4. $(\bigwedge_{\varphi \in \Psi} \varphi)^{\underline{M}} := \bigcap_{\varphi \in \Psi} \varphi^{\underline{M}}$, for $\Psi$ such that $\varnothing \neq \Psi \subseteq \Phi^\kappa$ and $|\Psi| < \kappa$.

So, defined as above, $\kappa$-expressions define measurable subsets of $M$. It is easy to check that $(\bigvee_{\varphi \in \Psi} \varphi)^{\underline{M}} := \bigcup_{\varphi \in \Psi} \varphi^{\underline{M}}$, for $\Psi$ such that $\varnothing \neq \Psi \subseteq \Phi^\kappa$ and $|\Psi| < \kappa$.

If no confusion may arise, we omit—with some abuse of notation—the superscript $\underline{M}$.

DEFINITION 12. For a $\kappa$-type space $\underline{M} := \langle M, \Sigma, (T_i)_{i \in I}, \theta \rangle$ on $S$ for player set $I$ and $m \in M$ define $D^\kappa(m)$, the $\kappa$-description of $m$, as

$$D^\kappa(m) := \{\varphi \in \Phi^\kappa | m \in \varphi^{\underline{M}}\}.$$

Again, we omit sometimes the superscript $\kappa$.

The next proposition says that type morphisms preserve $\kappa$-descriptions.

PROPOSITION 2. Let $\langle M, \Sigma, (T_i)_{i \in I}, \theta \rangle$ and $\langle N, \Sigma^N, (T_i^N)_{i \in I}, \theta^N \rangle$ be $\kappa$-type spaces on $S$ for player set $I$ and let $f : M \to N$ be a type morphism. Then, for all $m \in M$:

$$D(f(m)) = D(m).$$

PROOF. We show by induction on the formation of the expressions that $m \in \varphi^{\underline{M}}$ iff $f(m) \in \varphi^{\underline{N}}$:

(a) Let $E \in \Sigma_S$. We have $\theta^N(f(m)) = \theta(m)$, so $f(m) \in E^{\underline{N}}$ iff $m \in E^{\underline{M}}$.

(b) We have

$$f(m) \in (\neg\varphi)^{\underline{N}} \quad \text{iff} \quad f(m) \notin \varphi^{\underline{N}} \quad \text{iff} \quad m \notin \varphi^{\underline{M}} \quad \text{iff} \quad m \in (\neg\varphi)^{\underline{M}}.$$

(c) Let $\Psi$ be a nonempty set of expressions with $|\Psi| < \kappa$. Then

$$f(m) \in \left(\bigwedge_{\varphi \in \Psi} \varphi\right)^{\underline{N}} \quad \text{iff for all } \varphi \in \Psi : f(m) \in \varphi^{\underline{N}},$$

which is by the Induction hypothesis the case iff for all $\varphi \in \Psi : m \in \varphi^{\underline{M}}$, which is the case iff $m \in (\bigwedge_{\varphi \in \Psi} \varphi)^{\underline{M}}$.

(d) We have

$$f(m) \in (B_i^p(\varphi))^{\underline{N}} \quad \text{iff} \quad T_i^N(f(m))(\varphi^{\underline{N}}) \geq p \quad \text{iff} \quad T_i(m)(f^{-1}(\varphi^{\underline{N}})) \geq p.$$

By the Induction hypothesis: $f^{-1}(\varphi^{\underline{N}}) = \varphi^{\underline{M}}$. Hence $T_i(m)(f^{-1}(\varphi^{\underline{N}})) = T_i(m)(\varphi^{\underline{M}})$. We have $T_i(m)(\varphi^{\underline{M}}) \geq p$ iff $m \in (B_i^p(\varphi))^{\underline{M}}$. It follows that

$$f(m) \in (B_i^p(\varphi))^{\underline{N}} \quad \text{iff} \quad m \in (B_i^p(\varphi))^{\underline{M}}. \qquad \square$$

DEFINITION 13. Define $\Omega^\kappa$ to be the set of all $\kappa$-descriptions of states of the world in $\kappa$-type spaces on $S$ for player set $I$. For $\varphi \in \Phi^\kappa$ define

$$[\varphi] := \{\omega \in \Omega^\kappa | \varphi \in \omega\}.$$

Again, we omit sometimes the superscript $\kappa$.

REMARK 5. *The set of $\kappa$-descriptions $\Omega^\kappa$ is nonempty.*

PROOF. Let $M := \{m\}$ and choose $s \in S$. Set $\Sigma := \mathrm{Pow}(M)$, $T_i(m) := \delta_m$, for $i \in I$, and $\theta(m) := s$. Then

$$\langle M, \Sigma, (T_i)_{i \in I}, \theta \rangle$$

is a $\kappa$-type space (even a $*$-type space) on $S$ for player set $I$ and hence $D(m) \in \Omega^\kappa$. □

Obviously, we have $\Omega \setminus [\varphi] = [\neg\varphi]$ and $\bigcap_{\psi \in \Psi}[\psi] = [\bigwedge_{\psi \in \Psi} \psi]$, where $\varphi$ is a $\kappa$-expression and $\Psi$ is a nonempty set of $\kappa$-expressions with $|\Psi| < \kappa$. It follows that:

REMARK 6. *The set*

$$\Sigma_\Omega := \{[\varphi] | \varphi \in \Phi^\kappa\}$$

*is a $\kappa$-field on $\Omega$.*

LEMMA 4. *For every $\kappa$-type space $\underline{M}$ on $S$ for player set $I$ and for every $\varphi \in \Phi^\kappa$, the $\kappa$-description map $D \colon M \to \Omega$ satisfies*

$$D^{-1}([\varphi]) = \varphi\underline{^M}.$$

PROOF. Clear by the definition of $[\varphi]$. □

Note that Lemma 4 implies that $D$ is measurable.

PROPOSITION 3. *For every $i \in I$ there exists a function*

$$T_i^* \colon \Omega \to \Delta^\kappa(\Omega, \Sigma_\Omega)$$

*such that for every $\kappa$-type space $\underline{M}$ on $S$ for player set $I$ with $\kappa$-description map $D$ and every $m \in M$:*

$$T_i^*(D(m))(\cdot) = T_i(m)(D^{-1}(\cdot)).$$

PROOF. For $\omega \in \Omega$ choose a $\kappa$-type space $\underline{M}$ on $S$ for player set $I$ and $m \in M$ such that $D(m) = \omega$. For $[\varphi] \in \Sigma_\Omega$ define

$$T_i^*(\omega)([\varphi]) := T_i(m)(D^{-1}([\varphi])).$$

We have

$$T_i(m)(D^{-1}([\varphi])) = T_i(m)(\varphi^{\underline{M}}) = \sup\{p|B_i^p(\varphi) \in D(m)\},$$

so $T_i^*(\omega)([\varphi])$ depends just on $D(m)$ and is well-defined. By Remark 2, we have

$$T_i(m)(D^{-1}(\cdot)) \in \Delta^\kappa(\Omega, \Sigma_\Omega). \qquad \square$$

LEMMA 5. *There is a measurable function $\theta^* : \Omega \to S$ such that for every $\kappa$-type space $\underline{M}$ on $S$ for player set $I$ and every $m \in M$:*

$$\theta^*(D(m)) = \theta(m).$$

PROOF. Let

$$d_0(m) := \{E \in \Sigma_S | m \in \theta^{-1}(E)\}.$$

Obviously, $d_0(m) = D(m) \cap \Sigma_S$. By the properties of $(S, \Sigma_S)$, we have for all $s \in S : \{s\} = \bigcap_{s \in E \in \Sigma_S} E$. It follows for every $\kappa$-type space $\underline{M}'$ on $S$ for player set $I$ and $m' \in M'$ that

$$\theta(m') = s \quad \text{iff} \quad d_0(m') = \{E|s \in E\}.$$

For $\omega \in \Omega$ choose a $\kappa$-type space $\underline{M}$ on $S$ for player set $I$ and $m \in M$, such that $D(m) = \omega$. Define now $\theta^*(\omega) := \theta(m)$. Since $\theta(m)$ just depends on $D(m)$, $\theta^*(\omega)$ is well-defined.

It remains to show that $\theta^*$ is measurable: Let $E \in \Sigma_S$. We have

$$\theta^*(D(m)) \in E \quad \text{iff} \quad m \in \theta^{-1}(E) \quad \text{iff} \quad E \in D(m) \quad \text{iff} \quad D(m) \in [E].$$

It follows that $\theta^{*-1}(E) = [E]$. $\quad\square$

PROPOSITION 4.

$$\langle \Omega, \Sigma_\Omega, (T_i^*)_{i \in I}, \theta^* \rangle$$

*is a $\kappa$-type space on $S$ for player set $I$.*

PROOF. It remains to show:

1. For every $i \in I : T_i^*$ is measurable as a function from $\Omega$ to $\Delta^\kappa(\Omega, \Sigma_\Omega)$.
2. For every $i \in I$, $\omega \in \Omega$ and $A \in \Sigma_\Omega$: If

$$\{\omega' \in \Omega | T_i^*(\omega') = T_i^*(\omega)\} \subseteq A,$$

   then $T_i^*(\omega)(A) = 1$.

1. Since inverse images commute with unions, intersections and complements, it is enough to show that $T_i^{*-1}(b^p(E)) \in \Sigma_\Omega$, for

$$b^p(E) := \{\mu \in \Delta^\kappa(\Omega, \Sigma_\Omega) | \mu(E) \geq p\},$$

where $E \in \Sigma_\Omega$ and $p \in [0, 1]$. We have

$$T_i^{*-1}(b^p(E)) = \{\omega \in \Omega | T_i^*(\omega)(E) \geq p\}.$$

Since $E \in \Sigma_\Omega$, there is a $\kappa$-expression $\varphi$ such that $E = [\varphi]$. Note that if $p \in [0, 1]$ and $p = \sup\{q | B_i^q(\varphi) \in \omega\}$, then $B_i^p(\varphi) \in \omega$. This implies that

$$\omega \in T_i^{*-1}(b^p([\varphi])) \qquad \text{iff } B_i^p(\varphi) \in \omega.$$

It follows that $T_i^{*-1}(b^p(E)) = [B_i^p(\varphi)]$.

2. Let $\varphi$ be a $\kappa$-expression and

$$\{\omega' \in \Omega | T_i^*(\omega') = T_i^*(\omega)\} \subseteq [\varphi].$$

Choose a $\kappa$-type space $\underline{M}$ on $S$ for player set $I$ and $m \in M$ such that $D(m) = \omega$. Let $m' \in M$. If $T_i^*(D(m')) \neq T_i^*(D(m))$, then there is a $\kappa$-expression $\psi$ such that $T_i(m')(\psi^{\underline{M}}) \neq T_i(m)(\psi^{\underline{M}})$. It follows that

$$D(\{m' \in M | T_i(m') = T_i(m)\}) \subseteq \{\omega' \in \Omega | T_i^*(\omega') = T_i^*(\omega)\},$$

which implies

$$\{m' \in M | T_i(m') = T_i(m)\} \subseteq D^{-1}([\varphi]) = \varphi^{\underline{M}}.$$

So we have

$$1 = T_i(m)(\varphi^{\underline{M}}) = T_i(m)(D^{-1}([\varphi])) = T_i^*(\omega)([\varphi]). \qquad \square$$

Lemma 6. *The $\kappa$-description map*

$$D : \Omega \to \Omega$$

*is the identity.*

Proof. For $\omega \in \Omega$, we have

$$\omega = \{\varphi \in \Phi | \omega \in [\varphi]\}.$$

We have to show that for every $\kappa$-expression $\varphi$ and every $\omega \in \Omega : \omega \in \varphi^{\underline{\Omega}}$ iff $\omega \in [\varphi]$. We know this already if $\varphi = E$, where $E \in \Sigma_S$. It is obvious that $\Omega \setminus [\varphi] = [\neg\varphi]$, and that if $\Psi$ is a nonempty set of $\kappa$-expressions of cardinality $< \kappa$, then

$$\bigcap_{\varphi \in \Psi} [\varphi] = \left[\bigwedge_{\varphi \in \Psi} \varphi\right].$$

So it remains to show that $[\varphi] = \varphi^{\underline{\Omega}}$ implies $[B_i^p(\varphi)] = \overline{B}_i^p([\varphi])$. For $\omega \in \Omega$, choose a $\kappa$-type space $\underline{M}$ on $S$ for player set $I$ and $m \in M$ such that $D(m) = \omega$. We have

$$D(m) \in [B_i^p(\varphi)] \quad \text{iff} \quad B_i^p(\varphi) \in D(m) \quad \text{iff} \quad T_i(m)(\varphi^{\underline{M}}) \geq p.$$

But we have

$$T_i^*(\omega)([\varphi]) = T_i(m)D^{-1}([\varphi]) = T_i(m)(\varphi^{\underline{M}}).$$

This implies that $[B_i^p(\varphi)] = \overline{B}_i^p([\varphi])$. $\quad \square$

THEOREM 1. *The space*

$$\langle \Omega, \Sigma_\Omega, (T_i^*)_{i \in I}, \theta^* \rangle$$

*is a universal $\kappa$-type space on $S$ for player set $I$.*

PROOF. According to Lemma 4, for every $\kappa$-type space $\underline{M}$ on $S$ for player set $I$, the $\kappa$-description map $D: M \to \Omega$ is measurable, and according to Proposition 3 and Lemma 5, $D$ is a type morphism. It remains to show that it is the unique type morphism from $\underline{M}$ to $\underline{\Omega}$. But this is clear by Proposition 2 and Lemma 6. $\quad \square$

## 5. Spaces of arbitrary complexity.

Is there a cardinal $\kappa$, such that the $\kappa$-descriptions determine already the $\kappa'$-descriptions, for all cardinal numbers $\kappa' > \kappa$? In the sequel, using a probabilistic adaptation of the elegant "sober-drunk" example of Heifetz and Samet [9] (see that paper also for the "story" interpreting the mathematical structure), we construct, for every regular cardinal $\kappa'$, a $\kappa'$-type space (in fact even a $*$-type space), such that for every ordinal $\alpha < \kappa'$ there are at least two states of the world such that for every $\kappa'$-expression of depth $\leq \alpha$ this $\kappa'$-expression is true either in both states or in neither of the two, and yet there is a $\kappa'$-expression of depth $\alpha + 1$ that is true in one state and not in the other. Thus, choosing $\alpha$ and $\kappa'$ with $\kappa \leq \alpha < \alpha + 1 < \kappa'$, we answer the above question in the negative. Hence, it makes sense to consider $\kappa$-type spaces for every regular cardinality $\kappa$ whatsoever.

On top of that, this example will imply that, for at least two players and at least two states of nature, there is no universal $*$-type space and no universal $\infty$-type space (Theorem 4 and Corollary 1).

For this section, let $I := \{a, b\}$ be the set of players (the following analysis can be trivially extended to more than two players). We fix a set of states of nature $S = \{h, t\}$, consisting of the two possible outcomes of tossing a coin, $h$(ead) and $t$(ail).

To simplify the notation let us make the following

CONVENTION 2.    $\{i, j\} := \{a, b\}$, that is,

$$j = \begin{cases} a, & \text{if } i = b, \\ b, & \text{if } i = a. \end{cases}$$

The following three definitions and Definition 19 are taken from Heifetz and Samet [9].

DEFINITION 14.    Let $\alpha \geq 1$ be an ordinal. A *record of length* $\alpha$ is a sequence $r^\alpha = (r(\beta))_{\beta < \alpha}$ of numbers "0" and "1" such that for every limit ordinal $\lambda \leq \alpha$ there is an ordinal $\gamma < \lambda$ such that $r(\beta) = 0$ for all ordinals $\beta$ that satisfy $\gamma \leq \beta < \lambda$.

For every infinite ordinal $\gamma$ there are a unique natural number $n$ and a unique limit ordinal $\widehat{\lambda}$ such that $\gamma = \widehat{\lambda} + n$. We say $\gamma$ is *even* or *odd* according to whether $n$ is even or odd. If $\gamma$ is a finite ordinal, that is, a natural number, we take the usual notion of being even or odd.

DEFINITION 15.    Let $\alpha$ be an ordinal, $r^\alpha$ a record of length $\alpha$ and $\lambda$ a limit ordinal $\leq \alpha$. By the definition of a record, there is a minimal ordinal $o^\lambda(r^\alpha) < \lambda$ such that $r^\alpha(\beta) = 0$, for all $\beta$ with $o^\lambda(r^\alpha) \leq \beta < \lambda$.
Define $\lambda - \text{par}(r^\alpha)$, the *$\lambda$-parity* of $r^\alpha$, as

$$\lambda - \text{par}(r^\alpha) := \begin{cases} \text{even}, & \text{if } o^\lambda(r^\alpha) \text{ is even}, \\ \text{odd}, & \text{if } o^\lambda(r^\alpha) \text{ is odd}. \end{cases}$$

Note that by the definition of a record, $o^\lambda(r^\alpha)$ must be either 0 or a successor ordinal [i.e., $o^\lambda(r^\alpha) = \gamma + 1$, for some ordinal $\gamma$].

DEFINITION 16.    Let $\alpha$ be an ordinal. Define the spaces $W^\alpha$ by:

(i)  $W^0 := \{h, t\}$,
(ii) $W^\alpha := \{(w_0, w_a^\alpha, w_b^\alpha) | w_0 \in \{h, t\}, w_a^\alpha \text{ and } w_b^\alpha \text{ are records of length } \alpha\}$,
if $\alpha \geq 1$.

DEFINITION 17.    (i) If $0 < \beta \leq \alpha$ and $r^\alpha = (r(\xi))_{\xi < \alpha}$ is a record of length $\alpha$, then denote by $r^\alpha \lceil \beta$ the record $(r(\xi))_{\xi < \beta}$ of length $\beta$.
(ii) If $0 \leq \beta$ and $w^\alpha \in W^\alpha$, then define $w^\alpha \lceil 0 := w_0$.
(iii) If $0 < \beta \leq \alpha$ and $w^\alpha \in W^\alpha$, then define $w^\alpha \lceil \beta := (w_0, w_a^\alpha \lceil \beta, w_b^\alpha \lceil \beta)$.

By the definition, it is obvious that $w^\alpha \lceil \beta \in W^\beta$, for every $\beta < \alpha$.

DEFINITION 18.    Let $0 \leq \beta \leq \alpha$. Define

$$\pi_{\beta, \alpha} : W^\alpha \to W^\beta$$

by $\pi_{\beta, \alpha}(w^\alpha) := w^\alpha \lceil \beta$.

It is obvious that $\pi_{\xi,\beta}(\pi_{\beta,\alpha}(w^\alpha)) = \pi_{\xi,\alpha}(w^\alpha)$, for $0 \le \xi \le \beta \le \alpha$.

REMARK 7. *Let $0 \le \beta < \alpha$, $w^\beta \in W^\beta$, and $i \in \{a, b\}$. Then there are $w^\alpha$, $u^\alpha \in W^\alpha$ such that*

$$w^\alpha \lceil \beta = u^\alpha \lceil \beta = w^\beta \quad and \quad 0 = w_i^\alpha(\beta) \ne u_i^\alpha(\beta) = 1.$$

*In particular, it follows that $\pi_{\beta,\alpha} : W^\alpha \to W^\beta$ is onto.*

We define for each player $i$ a partition of the space $W^\alpha$. Two states are in the same element of player $i$'s partition if he cannot distinguish them. That is, he has the same information (the same beliefs) in both states. Let $w^\alpha = (w_0, w_a^\alpha, w_b^\alpha)$ be a state in $W^\alpha$. The element of $i$'s partition that contains this state is defined as follows:

DEFINITION 19. Let $\alpha$ be an ordinal $> 0$ and $w^\alpha \in W^\alpha$. We define:

$P_i(w^\alpha) := \{(v_0, v_a^\alpha, v_b^\alpha) \in W^\alpha | v_i^\alpha = w_i^\alpha,$

$\qquad\qquad w_i^\alpha(0) = 1$ implies $v_0 = w_0$,

$\qquad\qquad$ for all $\beta$ such that $\beta + 1 < \alpha :$

$\qquad\qquad w_i^\alpha(\beta + 1) = 1$ implies $v_j^\alpha(\beta) = w_j^\alpha(\beta)$,

$\qquad\qquad$ for every limit ordinal $\lambda < \alpha :$

$\qquad\qquad w_i^\alpha(\lambda) = 1$ implies $\lambda - \mathrm{par}(v_j^\alpha) = \lambda - \mathrm{par}(w_j^\alpha)\}$.

REMARK 8. (i) *Let $\alpha$ be an ordinal $> 0$. The set $\{P_i(w^\alpha) | w^\alpha \in W^\alpha\}$ is a partition of $W^\alpha$ and $w^\alpha \in P_i(w^\alpha)$.*

(ii) *Let $0 < \beta < \alpha$ and $u^\alpha \in P_i(w^\alpha)$. Then $u^\alpha \lceil \beta \in P_i(w^\alpha \lceil \beta)$, and hence $\pi_{\beta,\alpha}^{-1}(P_i(w^\alpha \lceil \beta)) \supseteq P_i(w^\alpha)$.*

It is easy to see that if $\alpha$ is an infinite ordinal number, then the cardinality of $W^\alpha$ is the same as the cardinality of $\alpha$. [To see that, in the case of an infinite $\alpha$, the cardinality of $W^\alpha$ does not exceed that of $\alpha$, note that the definition of a record implies that there are only finitely many $\beta < \alpha$ such that $r^\alpha(\beta) = 1$. Consider, assuming the contrary, the minimal $\gamma \le \alpha$ such that there are infinitely many $\beta < \gamma$ with $r^\alpha(\beta) = 1$.] Therefore, we have:

REMARK 9. *Let $\kappa$ be an infinite cardinal number. Then $|W^\kappa| = \kappa$.*

NOTATION 2. (a) For $0 \le \alpha$ and $w_0 \in \{h, t\}$, we denote by $[X_0^\alpha = w_0]$ the set $\{u^\alpha \in W^\alpha | u_0 = w_0\}$.

(b) For $\beta < \alpha$, $i \in \{a, b\}$ and $w_i^\alpha(\beta) \in \{0, 1\}$, we denote by $[X_i^\alpha(\beta) = w_i^\alpha(\beta)]$ the set $\{u^\alpha \in W^\alpha | u_i^\alpha(\beta) = w_i^\alpha(\beta)\}$.

(c) For a limit ordinal $\lambda \le \alpha$ and $\lambda - \mathrm{par}(w_i^\alpha) \in \{\text{even, odd}\}$, we denote by $[\lambda - \mathrm{par}(X_i^\alpha) = \lambda - \mathrm{par}(w_i^\alpha)]$ the set $\{u^\alpha \in W^\alpha | \lambda - \mathrm{par}(u_i^\alpha) = \lambda - \mathrm{par}(w_i^\alpha)\}$.

REMARK 10. *Let $0 \leq \alpha \leq \gamma$ and $w_0 \in \{h, t\}$. Then:*

(i) $\pi_{\alpha,\gamma}^{-1}([X_0^\alpha = w_0]) = [X_0^\gamma = w_0]$.

(ii) *If $\beta < \alpha$, $i \in \{a, b\}$ and $w^\gamma \in W^\gamma$, then*

$$w_i^\gamma(\beta) = (w_i^\gamma \lceil \alpha)(\beta)$$

*and*

$$\pi_{\alpha,\gamma}^{-1}([X_i^\alpha(\beta) = (w_i^\gamma \lceil \alpha)(\beta)]) = [X_i^\gamma(\beta) = w_i^\gamma(\beta)].$$

(iii) *If $i \in \{a, b\}$, $w^\gamma \in W^\gamma$, and if $\lambda$ is a limit ordinal such that $\lambda \leq \alpha$, then*

$$\lambda - \mathrm{par}(w_i^\gamma) = \lambda - \mathrm{par}(w_i^\gamma \lceil \alpha)$$

*and*

$$\pi_{\alpha,\gamma}^{-1}([\lambda - \mathrm{par}(X_i^\alpha) = \lambda - \mathrm{par}(w_i^\gamma \lceil \alpha)]) = [\lambda - \mathrm{par}(X_i^\gamma) = \lambda - \mathrm{par}(w_i^\gamma)].$$

Now, let $\kappa$ be a fixed regular cardinal.

THEOREM 2. *For $i \in \{a, b\}$, there is a function $T_i \colon W^\kappa \longrightarrow \Delta(W^\kappa, \mathrm{Pow}(W^\kappa))$ with the following properties:*

(a) $T_i(u^\kappa) = T_i(w^\kappa)$, *for $u^\kappa \in P_i(w^\kappa)$,*

(b) $T_i(w^\kappa)(P_i(w^\kappa)) = 1$,

(c) $T_i(w^\kappa)([X_0^\kappa = w_0]) = \begin{cases} 1, & \text{if } w_i^\kappa(0) = 1, \\ \frac{1}{2}, & \text{if } w_i^\kappa(0) = 0, \end{cases}$

(d) *for $\beta < \kappa$:*

$$T_i(w^\kappa)([X_j^\kappa(\beta) = w_j^\kappa(\beta)]) = \begin{cases} 1, & \text{if } w_i^\kappa(\beta + 1) = 1, \\ \frac{1}{2}, & \text{if } w_i^\kappa(\beta + 1) = 0, \end{cases}$$

(e) *for $\lambda < \kappa$ such that $\lambda$ is a limit ordinal:*

$$T_i(w^\kappa)([\lambda - \mathrm{par}(X_j^\kappa) = \lambda - \mathrm{par}(w_j^\kappa)]) = \begin{cases} 1, & \text{if } w_i^\kappa(\lambda) = 1, \\ \frac{1}{2}, & \text{if } w_i^\kappa(\lambda) = 0, \end{cases}$$

(f) *for $0 \leq \beta < \alpha < \kappa$ and $u^\kappa, w^\kappa \in W^\kappa$: if $E^\beta \subseteq W^\beta$ and $u^\kappa \lceil \alpha = w^\kappa \lceil \alpha$, then*

$$T_i(u^\kappa)(\pi_{\beta,\kappa}^{-1}(E^\beta)) = T_i(w^\kappa)(\pi_{\beta,\kappa}^{-1}(E^\beta)).$$

Theorem 2 will be proved in Section 6.

REMARK 11. *For $w^\kappa \in W^\kappa$, define*

$$\theta(w^\kappa) := w_0.$$

*From the first two points of Theorem 2 and the fact that $T_i(w^\kappa)$ is a finitely additive probability measure defined on $(W^\kappa, \mathrm{Pow}(W^\kappa))$, it follows that*

$$\langle W^\kappa, (T_i)_{i \in \{a,b\}}, \theta \rangle$$

*is a $*$-type space on $S = \{h, t\}$ for player set $\{a, b\}$.*

Next, by induction on the formation of $\kappa$-expressions, we define the *depth* of a $\kappa$-expression. The depth of a $\kappa$-expression is an ordinal number that measures how complex that expression is with respect to the players' beliefs operators.

DEFINITION 20. (i) If $E \in \Sigma_S$, then $\mathrm{dp}(E) := 0$,

(ii) if $0 \leq p \leq 1$, $i \in I$ and if $\varphi$ is a $\kappa$-expression, then

$$\mathrm{dp}(B_i^p(\varphi)) := \mathrm{dp}(\varphi) + 1,$$

(iii) if $\varphi$ is a $\kappa$-expression, then $\mathrm{dp}(\neg \varphi) := \mathrm{dp}(\varphi)$,

(iv) if $\Psi$ is a set of $\kappa$-expressions such that $|\Psi| < \kappa$, then

$$\mathrm{dp}\left( \bigwedge_{\varphi \in \Psi} \varphi \right) := \sup\{\mathrm{dp}(\varphi) | \varphi \in \Psi\}.$$

It is easy to see that, since $\kappa$ is regular, the depth of a $\kappa$-expression is always strictly smaller than $\kappa$.

The following Lemmas 7–9 will be proved in the Appendix.

LEMMA 7. *Let $\alpha \leq \kappa$, $w^\kappa, u^\kappa \in W^\kappa$ and $w^\kappa \lceil \alpha = u^\kappa \lceil \alpha$.*
*Then, for all $\kappa$-expressions $\varphi$ such that $\mathrm{dp}(\varphi) \leq \alpha$:*

$$w^\kappa \in \varphi^{W^\kappa} \quad iff \quad u^\kappa \in \varphi^{W^\kappa}.$$

LEMMA 8. *In the $*$-type space $\langle W^\kappa, (T_i)_{i \in \{a,b\}}, \theta \rangle$, we have:*

$$[X_i^\kappa(0) = 1] = \overline{B}_i^1([X_0^\kappa = h]) \cup \overline{B}_i^1([X_0^\kappa = t]),$$

$$[X_i^\kappa(\beta + 1) = 1] = \overline{B}_i^1([X_j^\kappa(\beta) = 1]) \cup \overline{B}_i^1([X_j^\kappa(\beta) = 0])$$

$$\text{for all ordinals } \beta < \kappa,$$

$$[X_i^\kappa(\lambda) = 1] = \overline{B}_i^1([\lambda - \mathrm{par}(X_j^\kappa) = \mathrm{even}]) \cup \overline{B}_i^1([\lambda - \mathrm{par}(X_j^\kappa) = \mathrm{odd}])$$

$$\text{for all limit ordinals } \lambda < \kappa.$$

LEMMA 9. *In the $*$-type space $\langle W^\kappa, (T_i)_{i \in \{a,b\}}, \theta \rangle$, we have:*

(i)

$$\{h\}^{W^\kappa} = [X_0^\kappa = h],$$

$$\{t\}^{W^\kappa} = [X_0^\kappa = t],$$

$$\mathrm{dp}(\{h\}) = \mathrm{dp}(\{t\}) = 0.$$

(ii) *For every* $i \in \{a, b\}$ *and* $\beta$ *such that* $0 \leq \beta < \kappa$, *there are* $\kappa$-*expressions* $\varphi_i^0(\beta)$ *and* $\varphi_i^1(\beta)$ *with*

$$\mathrm{dp}(\varphi_i^0(\beta)) = \mathrm{dp}(\varphi_i^1(\beta)) = \beta + 1$$

*such that*

$$(\varphi_i^1(\beta))^{W^\kappa} = [X_i^\kappa(\beta) = 1]$$

*and*

$$(\varphi_i^0(\beta))^{W^\kappa} = [X_i^\kappa(\beta) = 0].$$

(iii) *For every* $i \in \{a, b\}$ *and limit ordinal* $\lambda < \kappa$, *there are* $\kappa$-*expressions* $\varphi_i^{\mathrm{even}}(\lambda)$ *and* $\varphi_i^{\mathrm{odd}}(\lambda)$ *with*

$$\mathrm{dp}(\varphi_i^{\mathrm{even}}(\lambda)) = \mathrm{dp}(\varphi_i^{\mathrm{odd}}(\lambda)) = \lambda$$

*such that*

$$(\varphi_i^{\mathrm{even}}(\lambda))^{W^\kappa} = [\lambda - \mathrm{par}(X_i^\kappa) = \mathrm{even}]$$

*and*

$$(\varphi_i^{\mathrm{odd}}(\lambda))^{W^\kappa} = [\lambda - \mathrm{par}(X_i^\kappa) = \mathrm{odd}].$$

THEOREM 3.   *For every ordinal* $\alpha < \kappa$ *there are* $u^\kappa, w^\kappa \in W^\kappa$ *such that:*

1. *For all* $\kappa$-*expressions* $\varphi$ *with* $\mathrm{dp}(\varphi) \leq \alpha$:

$$u^\kappa \in \varphi^{W^\kappa} \quad iff \quad w^\kappa \in \varphi^{W^\kappa}.$$

2. *There is a* $\kappa$-*expression* $\psi$ *with* $\mathrm{dp}(\psi) = \alpha + 1$ *such that*

$$u^\kappa \in \psi^{W^\kappa} \quad and \quad w^\kappa \in (\neg \psi)^{W^\kappa}.$$

PROOF.   Let $\alpha < \kappa$ and $i \in \{a, b\}$. By the definition of $W^\kappa$, there are $u^\kappa$, $w^\kappa \in W^\kappa$ such that $u^\kappa {\restriction} \alpha = w^\kappa {\restriction} \alpha$ and $1 = u_i^\kappa(\alpha) \neq w_i^\kappa(\alpha) = 0$. The first point follows now by Lemma 7. By Lemma 9, it follows that $u^\kappa \in \varphi_i^1(\alpha)^{W^\kappa}$, $w^\kappa \in (\neg\varphi_i^1(\alpha))^{W^\kappa}$, and $\mathrm{dp}(\varphi_i^1(\alpha)) = \alpha + 1$.   □

Note that Lemma 7 and the proof of Theorem 3 show that $u_0$ and the levels up to and including level $\alpha$ of $u_i^\kappa$ and of $u_j^\kappa$ determine which of the $\kappa$-expressions of depth $\alpha + 1$ belong to the $\kappa$-description of $u^\kappa$.

THEOREM 4. *Let $|I| \geq 2$ and $|S| \geq 2$. Then, there is no weak-universal $\infty$-type space on $S$ for player set $I$ and there is no weak-universal $*$-type space on $S$ for player set $I$.*

PROOF. Assume there is a weak-universal $\infty$-type space (a weak-universal $*$-type space, resp.)

$$\underline{U} = \langle U, \Sigma, (T_i^U)_{i \in I}, \theta^U \rangle$$

on $S$ for player set $I$. Then, the underlying set $U$ has a cardinality $|U|$. There is a regular cardinal number $\kappa > |U|$.

$$\underline{W}^\kappa = \langle W^\kappa, (T_i)_{i \in \{a,b\}}, \theta \rangle$$

is a $*$-type space on $\{h,t\}$ (and therefore a $\infty$-type space). Since $|\{h,t\}| = 2$, we can assume without loss of generality that $\{h,t\} \subseteq S$, and since $\Sigma_S = \mathrm{Pow}(S)$, that $\Sigma_S \supseteq \mathrm{Pow}(\{h,t\})$. Also, since $|I| \geq 2$, we can assume without loss of generality that $\{a,b\} \subseteq I$. For $i \in I \setminus \{a,b\}$ define $T_i(w^\kappa) := \delta_{w^\kappa}$, and view $\theta$ as a function from $W^\kappa$ to $S$. Then

$$\underline{W}_I^\kappa := \langle W^\kappa, (T_i)_{i \in I}, \theta \rangle$$

is a $*$-type space on $S$ (with player set $I$). Since every $*$-type space is a $\infty$-type space, $\underline{W}_I^\kappa$ is also a $\infty$-type space. According to the assumption, there is a type morphism $f : W^\kappa \to U$. Since both spaces are in particular $\kappa$-type spaces, this morphism preserves $\kappa$-descriptions. If $\varphi$ is a $\kappa$-expression in the "language" corresponding to the set of states of nature $\{h,t\}$ and the player set $\{a,b\}$, then $\varphi$ is also a $\kappa$-expression in the "language" corresponding to the set of states of nature $S$ and the player set $I$, and it is easy to check that for $w^\kappa \in W^\kappa$ we have $w^\kappa \in \varphi^{\underline{W}_I^\kappa}$ iff $w^\kappa \in \varphi^{\underline{W}^\kappa}$. So, by Lemma 9, it is still the case that two different states of $\underline{W}_I^\kappa$ have different $\kappa$-descriptions. Hence, since by Proposition 2, $f$ preserves $\kappa$-descriptions, $f$ is one-to-one. It follows that $|U| \geq |W^\kappa| = \kappa$, which is a contradiction to $|U| < \kappa$. □

COROLLARY 1. *Let $|I| \geq 2$ and $|S| \geq 2$. Then there is no universal $\infty$-type space on $S$ for player set $I$ and there is no universal $*$-type space on $S$ for player set $I$.*

**6. Construction of the $T_i$'s.** This section is devoted to the proof of Theorem 2, that is, the construction of the $T_i$'s mentioned there. Lemmas 10 and 12–17 needed for this construction are proved at the end of this section.

The construction will not be carried out at once. By a transfinite induction on $1 \leq \alpha \leq \kappa$, we endow $W^\alpha$ with fields $\mathcal{F}(i, w^\alpha)$ and finitely additive probability measures $T_i^\alpha(w^\alpha)$ on $(W^\alpha, \mathcal{F}(i, w^\alpha))$ such that the following Induction hypothesis is satisfied:

INDUCTION HYPOTHESIS (for $\alpha$).

1. $\mathcal{F}(i, w^\alpha) := [\bigcup_{\beta < \alpha} (\pi_{\beta,\alpha}^{-1}(\text{Pow}(W^\beta))), P_i(w^\alpha)]$.
2. For every ordinal $\beta$ with $1 \leq \beta < \alpha$ and every $E^\beta \in \mathcal{F}(i, w^\alpha \restriction \beta)$:

$$T_i^\alpha(w^\alpha)(\pi_{\beta,\alpha}^{-1}(E^\beta)) = T_i^\beta(w^\alpha \restriction \beta)(E^\beta),$$

   that is, $\text{marg}_{(W^\beta, \mathcal{F}(i, w^\alpha \restriction \beta))} T_i^\alpha(w^\alpha) = T_i^\beta(w^\alpha \restriction \beta)$.
3. $T_i^\alpha(w^\alpha) = T_i^\alpha(u^\alpha)$ for $u^\alpha \in P_i(w^\alpha)$.
4. $T_i^\alpha(w^\alpha)(P_i(w^\alpha)) = 1$.
5. $T_i^\alpha(w^\alpha)([X_0^\alpha = w_0]) = \begin{cases} 1, & \text{if } w_i^\alpha(0) = 1, \\ \frac{1}{2}, & \text{if } w_i^\alpha(0) = 0. \end{cases}$
6. For $\beta$ such that $\beta + 1 < \alpha$:

$$T_i^\alpha(w^\alpha)([X_j^\alpha(\beta) = w_j^\alpha(\beta)]) = \begin{cases} 1, & \text{if } w_i^\alpha(\beta + 1) = 1, \\ \frac{1}{2}, & \text{if } w_i^\alpha(\beta + 1) = 0. \end{cases}$$

7. For $\lambda < \alpha$ such that $\lambda$ is a limit ordinal:

$$T_i^\alpha(w^\alpha)([\lambda - \text{par}(X_j^\alpha) = \lambda - \text{par}(w_j^\alpha)]) = \begin{cases} 1, & \text{if } w_i^\alpha(\lambda) = 1, \\ \frac{1}{2}, & \text{if } w_i^\alpha(\lambda) = 0. \end{cases}$$

Let $1 \leq \gamma \leq \alpha$. Then $\bigcup_{\beta < \gamma} \pi_{\beta,\alpha}^{-1}(\text{Pow}(W^\beta))$ is a field on $W^\alpha$. Hence, by definition, $\mathcal{F}(i, w^\alpha)$ is also a field on $W^\alpha$. Since inverse images commute with complements, arbitrary unions and intersections, we have:

REMARK 12. *Let $1 \leq \gamma \leq \alpha$. Then,*

$$\pi_{\gamma,\alpha}^{-1}(\mathcal{F}(i, w^\gamma)) = \left[ \bigcup_{\beta < \gamma} \pi_{\beta,\alpha}^{-1}(\text{Pow}(W^\beta)), \pi_{\gamma,\alpha}^{-1}(P_i(w^\gamma)) \right]$$

*is a field on $W^\alpha$. In particular,*

$$\pi_{\beta+1,\alpha}^{-1}(\mathcal{F}(i, w^\alpha \restriction \beta + 1)) = [\pi_{\beta,\alpha}^{-1}(\text{Pow}(W^\beta)), \pi_{\beta+1,\alpha}^{-1}(P_i(w^\alpha \restriction \beta + 1))]$$

*is a field on $W^\alpha$, for $\beta + 1 \leq \alpha$. Note also that $\pi_{\gamma,\alpha}^{-1}(\mathcal{F}(i, w^\alpha \restriction \gamma)) \subseteq \mathcal{F}(i, w^\alpha)$, for $1 \leq \gamma \leq \alpha$.*

REMARK 13. *Let $1 \leq \beta \leq \alpha \leq \kappa$ and let $u^\alpha, w^\alpha \in W^\alpha$ such that $u^\alpha \in P_i(w^\alpha)$ [and hence $P_i(u^\alpha) = P_i(w^\alpha)$]. Then:*

(a) *We have*

$$\mathcal{F}(i, w^\alpha) = \mathcal{F}(i, u^\alpha),$$

$$\mathcal{F}(i, w^\alpha \restriction \beta) = \mathcal{F}(i, u^\alpha \restriction \beta).$$

(b) *If*

$$T_i^\alpha(w^\alpha) = T_i^\alpha(u^\alpha),$$

$$T_i^\beta(w^\alpha \restriction \beta) = T_i^\beta(u^\alpha \restriction \beta),$$

$$\mathrm{marg}_{(W^\beta, \mathcal{F}(i, w^\alpha \restriction \beta))} T_i^\alpha(w^\alpha) = T_i^\beta(w^\alpha \restriction \beta),$$

*then we obviously also have*

$$\mathrm{marg}_{(W^\beta, \mathcal{F}(i, u^\alpha \restriction \beta))} T_i^\alpha(u^\alpha) = T_i^\beta(u^\alpha \restriction \beta).$$

Before we begin the construction of the $T_i$'s, we have to provide some lemmas that are needed to carry out this construction. These lemmas will guarantee that the induction can be done maintaining the conditions 1–7 of the Induction hypothesis. They will be proved at the end of this section.

LEMMA 10. *Let $\gamma$ be an ordinal $> 0$, $\alpha = \gamma + 1$, $w^\alpha \in W^\alpha$, and*

$$E \in \left[ \left[ \bigcup_{\beta < \gamma} \pi_{\beta,\alpha}^{-1}(\mathrm{Pow}(W^\beta)), \pi_{\gamma,\alpha}^{-1}(P_i(w^\alpha \restriction \gamma)) \right], P_i(w^\alpha) \right].$$

*Then there are a $\beta < \gamma$ and $A_\beta, C_\beta, D_\beta \in \mathrm{Pow}(W^\beta)$ such that*

$$E = (\pi_{\beta,\alpha}^{-1}(A_\beta) \cap P_i(w^\alpha)) \cup (\pi_{\beta,\alpha}^{-1}(C_\beta) \cap \pi_{\gamma,\alpha}^{-1}(P_i(w^\alpha \restriction \gamma)) \cap (W^\alpha \setminus P_i(w^\alpha)))$$

$$\cup (\pi_{\beta,\alpha}^{-1}(D_\beta) \cap (W^\alpha \setminus \pi_{\gamma,\alpha}^{-1}(P_i(w^\alpha \restriction \gamma)))).$$

For further reference, we cite here (in a slightly changed formulation and in our notation) Lemma 3.2 of [9]:

LEMMA 11. *Let $v^\alpha, w^\alpha \in W^\alpha$, where $v^\alpha \restriction \gamma + 1 \in P_i(w^\alpha \restriction \gamma + 1)$, for some $\gamma < \alpha$. Then there is a $u^\alpha \in P_i(w^\alpha)$ such that $u^\alpha \restriction \gamma = v^\alpha \restriction \gamma$.*

LEMMA 12. *Let $\lambda$ be a limit ordinal, $\alpha = \lambda + 1$, $w^\alpha \in W^\alpha$, $w_i^\alpha(\lambda) = 0$ and $E = \pi_{\beta,\alpha}^{-1}(E_\beta)$, where $E_\beta \subseteq W^\beta$ for a $\beta < \lambda$. Then:*

(a) *If $v^\alpha \in E \cap P_i(w^\alpha)$, then there is a $u^\alpha \in E \cap P_i(w^\alpha)$ such that $\lambda - \mathrm{par}(u_j^\alpha) \neq \lambda - \mathrm{par}(v_j^\alpha)$.*

(b) *If $v^\alpha \in E \cap \pi_{\lambda,\alpha}^{-1}(P_i(w^\alpha \restriction \lambda)) \cap (W^\alpha \setminus P_i(w^\alpha))$, then there is a $u^\alpha \in E \cap \pi_{\lambda,\alpha}^{-1}(P_i(w^\alpha \restriction \lambda)) \cap (W^\alpha \setminus P_i(w^\alpha))$ such that $\lambda - \mathrm{par}(u_j^\alpha) \neq \lambda - \mathrm{par}(v_j^\alpha)$.*

(c) *If $v^\alpha \in E \cap (W^\alpha \setminus \pi_{\lambda,\alpha}^{-1}(P_i(w^\alpha \restriction \lambda)))$, then there is a $u^\alpha \in E \cap (W^\alpha \setminus \pi_{\lambda,\alpha}^{-1}(P_i(w^\alpha \restriction \lambda)))$ such that $\lambda - \mathrm{par}(u_j^\alpha) \neq \lambda - \mathrm{par}(v_j^\alpha)$.*

LEMMA 13. *Let $\beta$ be an ordinal, $\alpha = (\beta + 1) + 1$, $w^\alpha \in W^\alpha$, $w_i^\alpha(\beta + 1) = 0$ and $E = \pi_{\beta,\alpha}^{-1}(E_\beta)$, such that $E_\beta \subseteq W^\beta$. Then:*

(a) *If $v^\alpha \in E \cap P_i(w^\alpha)$, then there is a $u^\alpha \in E \cap P_i(w^\alpha)$ such that $u_j^\alpha(\beta) \neq v_j^\alpha(\beta)$.*

(b) *If $v^\alpha \in E \cap \pi_{\beta+1,\alpha}^{-1}(P_i(w^\alpha \lceil \beta + 1)) \cap (W^\alpha \setminus P_i(w^\alpha))$, then there is a $u^\alpha \in E \cap \pi_{\beta+1,\alpha}^{-1}(P_i(w^\alpha \lceil \beta + 1)) \cap (W^\alpha \setminus P_i(w^\alpha))$ such that $u_j^\alpha(\beta) \neq v_j^\alpha(\beta)$.*

(c) *If $v^\alpha \in E \cap (W^\alpha \setminus \pi_{\beta+1,\alpha}^{-1}(P_i(w^\alpha \lceil \beta + 1)))$, then there is a $u^\alpha \in E \cap (W^\alpha \setminus \pi_{\beta+1,\alpha}^{-1}(P_i(w^\alpha \lceil \beta + 1)))$ such that $u_j^\alpha(\beta) \neq v_j^\alpha(\beta)$.*

LEMMA 14.  *Let $\lambda$ be a limit ordinal, $\alpha = \lambda + 1$, $w^\alpha \in W^\alpha$, $w_i^\alpha(\lambda) = 0$ and*

$$E \in \left[ \left[ \bigcup_{\beta < \lambda} \pi_{\beta,\alpha}^{-1}(\mathrm{Pow}(W^\beta)), \pi_{\lambda,\alpha}^{-1}(P_i(w^\alpha \lceil \lambda)) \right], P_i(w^\alpha) \right].$$

*Then:*

(a) *If $v^\alpha \in E$, then there is a $u^\alpha \in E$ such that $\lambda - \mathrm{par}(u_j^\alpha) \neq \lambda - \mathrm{par}(v_j^\alpha)$.*

(b) *If $v^\alpha \in W^\alpha \setminus E$, then there is a $u^\alpha \in W^\alpha \setminus E$ such that $\lambda - \mathrm{par}(u_j^\alpha) \neq \lambda - \mathrm{par}(v_j^\alpha)$.*

(c) *If $E \supseteq [\lambda - \mathrm{par}(X_j^\alpha) = \lambda - \mathrm{par}(w_j^\alpha)]$, then $E = W^\alpha$.*

(d) *If $E \subseteq [\lambda - \mathrm{par}(X_j^\alpha) = \lambda - \mathrm{par}(w_j^\alpha)]$, then $E = \varnothing$.*

LEMMA 15.  *Let $\beta$ be an ordinal, $\alpha = (\beta + 1) + 1$, $w^\alpha \in W^\alpha$, $w_i^\alpha(\beta + 1) = 0$ and*

$$E \in [[\pi_{\beta,\alpha}^{-1}(\mathrm{Pow}(W^\beta)), \pi_{\beta+1,\alpha}^{-1}(P_i(w^\alpha \lceil \beta + 1))], P_i(w^\alpha)].$$

*Then:*

(a) *If $v^\alpha \in E$, then there is a $u^\alpha \in E$ such that $u_j^\alpha(\beta) \neq v_j^\alpha(\beta)$.*

(b) *If $v^\alpha \in W^\alpha \setminus E$, then there is a $u^\alpha \in W^\alpha \setminus E$ such that $u_j^\alpha(\beta) \neq v_j^\alpha(\beta)$.*

(c) *If $E \supseteq [X_j^\alpha(\beta) = w_j^\alpha(\beta)]$, then $E = W^\alpha$.*

(d) *If $E \subseteq [X_j^\alpha(\beta) = w_j^\alpha(\beta)]$, then $E = \varnothing$.*

LEMMA 16.  *Let $\gamma$ be an ordinal $> 0$, $\alpha = \gamma + 1$, $w^\alpha \in W^\alpha$ and*

$$E \in \left[ \bigcup_{\beta < \gamma} \pi_{\beta,\alpha}^{-1}(\mathrm{Pow}(W^\beta)), \pi_{\gamma,\alpha}^{-1}(P_i(w^\alpha \lceil \gamma)) \right]$$

*such that $E \supseteq P_i(w^\alpha)$. Then*

$$E \supseteq \pi_{\gamma,\alpha}^{-1}(P_i(w^\alpha \lceil \gamma)).$$

LEMMA 17.  *Let $\lambda$ be a limit ordinal, $w^\lambda \in W^\lambda$, $\beta < \lambda$ and $E^\beta \subseteq W^\beta$ such that $\pi_{\beta,\lambda}^{-1}(E^\beta) \supseteq P_i(w^\lambda)$. Then $\pi_{\beta,\lambda}^{-1}(E^\beta) \supseteq \pi_{\beta+1,\lambda}^{-1}(P_i(w^\lambda \lceil \beta + 1))$.*

PROOF OF THEOREM 2.  We construct the $T_i$'s by transfinite induction on $1 \leq \alpha \leq \kappa$, such that the conditions 1–7 of the Induction hypothesis at the beginning of this section are satisfied:

*Step* $\alpha = 1$. We have $W^0 = \{h, t\}$. Let $w^1 \in W^1$. Define

$$T_i^{<1}(w^1)([X_0^1 = w_0]) := \begin{cases} 1, & \text{if } w_i^1(0) = 1, \\ \frac{1}{2}, & \text{if } w_i^1(0) = 0, \end{cases}$$

$$T_i^{<1}(w^1)([X_0^1 \neq w_0]) := 1 - T_i^{<1}(w^1)([X_0^1 = w_0]),$$

$$T_i^{<1}(w^1)(\varnothing) := 0,$$
$$T_i^{<1}(w^1)(W^1) := 1.$$

It is clear that by this definition, $T_i^{<1}(w^1)$ is a probability measure on

$$(W^1, \pi_{0,1}^{-1}(\mathrm{Pow}(W^0))).$$

Let $E^0 \subseteq W^0$ such that $P_i(w^1) \subseteq \pi_{0,1}^{-1}(E^0)$.

1. Case: $w_i^1(0) = 1$. Then $E^0 = W^0$ or $E^0 = [X_0^1 = w_0]$, hence the outer measure $T_i^{<1}(w^1)^*(P_i(w^1))$ is equal to 1.
2. Case: $w_i^1(0) = 0$. Then $E^0 = W^0$ and the outer measure $T_i^{<1}(w^1)^*(P_i(w^1))$ is equal to 1.

For $u^1 \in P_i(w^1)$, we have in both cases that $P_i(u^1) = P_i(w^1)$ and $T_i^{<1}(u^1) = T_i^{<1}(w^1)$. For each $P_i(u^1)$ such that $u^1 \in W^1$, choose a representing element $w^1 \in P_i(u^1) = P_i(w^1)$. By the Łoś–Marczewski theorem, we can extend $T_i^{<1}(w^1)$ to a finitely additive probability measure $T_i^1(w^1)$ on the field

$$\mathcal{F}(i, w^1) = [\pi_{0,1}^{-1}(\mathrm{Pow}(W^0)), P_i(w^1)]$$

such that $T_i^1(w^1)(P_i(w^1)) = 1$. Define $T_i^1(u^1) := T_i^1(w^1)$, for all $u^1 \in P_i(w^1)$.

Note that $\mathcal{F}(i, u^1) = \mathcal{F}(i, w^1)$, for $u^1 \in P_i(w^1)$. $T_i^1(u^1)$ and $\mathcal{F}(i, u^1)$ satisfy the conditions 1–7 of the Induction hypothesis.

*Step* $\alpha = (\beta + 1) + 1$, *for* $0 \leq \beta < \kappa$. For each $P_i(u^\alpha)$ such that $u^\alpha \in W^\alpha$, choose a representing element $w^\alpha \in P_i(u^\alpha) = P_i(w^\alpha)$.

Let $T_i^{<\alpha}(w^\alpha)$ be the finitely additive probability measure defined on the field

$$\pi_{\beta+1,\alpha}^{-1}(\mathcal{F}(i, w^\alpha \lceil \beta + 1)) = [\pi_{\beta,\alpha}^{-1}(\mathrm{Pow}(W^\beta)), \pi_{\beta+1,\alpha}^{-1}(P_i(w^\alpha \lceil \beta + 1))],$$

which is induced by $T_i^{\beta+1}(w^\alpha \lceil \beta + 1)$ (as defined in Lemma 1). According to Lemma 16 and the Induction hypothesis, we have for the outer measure of $P_i(w^\alpha): T_i^{<\alpha}(w^\alpha)^*(P_i(w^\alpha)) = 1$. So, by the Łoś–Marczewski theorem, we can extend $T_i^{<\alpha}(w^\alpha)$ to a finitely additive probability measure $\widetilde{T}_i^\alpha(w^\alpha)$ defined on the field

$$\widetilde{\mathcal{F}}(i, w^\alpha) := [[\pi_{\beta,\alpha}^{-1}(\mathrm{Pow}(W^\beta)), \pi_{\beta+1,\alpha}^{-1}(P_i(w^\alpha \lceil \beta + 1))], P_i(w^\alpha)]$$

such that $\widetilde{T}_i^\alpha(w^\alpha)(P_i(w^\alpha)) = 1$.

1. Case: $w_i^\alpha(\beta+1) = 1$. Then

$$\pi_{\beta+1,\alpha}^{-1}([X_j^{\beta+1}(\beta) = (w_j^\alpha\lceil\beta+1)(\beta)]) = [X_j^\alpha(\beta) = w_j^\alpha(\beta)] \supseteq P_i(w^\alpha).$$

   By Lemma 3, extend $\widetilde{T}_i^\alpha(w^\alpha)$ to a finitely additive probability measure $T_i^\alpha(w^\alpha)$ on the field

$$[\pi_{\beta+1,\alpha}^{-1}(\mathrm{Pow}(W^{\beta+1})), P_i(w^\alpha)] = \mathcal{F}(i, w^\alpha).$$

   By the above, we have

$$T_i^\alpha(w^\alpha)([X_j^\alpha(\beta) = w_j^\alpha(\beta)]) = 1.$$

   Define now $T_i^\alpha(u^\alpha) := T_i^\alpha(w^\alpha)$, for all $u^\alpha \in P_i(w^\alpha)$. Note that, for $u^\alpha \in P_i(w^\alpha)$, we have $\widetilde{\mathcal{F}}(i, u^\alpha) = \widetilde{\mathcal{F}}(i, w^\alpha)$, $\mathcal{F}(i, u^\alpha) = \mathcal{F}(i, w^\alpha)$ and $u_i^\alpha(\beta+1) = w_i^\alpha(\beta+1) = 1$ and hence, $u_j^\alpha(\beta) = w_j^\alpha(\beta)$. It is now easy to check that $T_i^\alpha(u^\alpha)$ and $\mathcal{F}(i, u^\alpha)$ satisfy the conditions 1–7 of the Induction hypothesis.

2. Case: $w_i^\alpha(\beta+1) = 0$. By Lemma 15 and the Induction hypothesis, we have for the outer measure of $[X_j^\alpha(\beta) = w_j^\alpha(\beta)]$:

$$\widetilde{T}_i^\alpha(w^\alpha)^*([X_j^\alpha(\beta) = w_j^\alpha(\beta)]) = 1,$$

   and for the inner measure

$$\widetilde{T}_i^\alpha(w^\alpha)_*([X_j^\alpha(\beta) = w_j^\alpha(\beta)]) = 0.$$

   By the Łoś–Marczewski theorem, we can extend $\widetilde{T}_i^\alpha(w^\alpha)$ to a finitely additive probability measure $\widehat{T}_i^\alpha(w^\alpha)$ on the field

$$[\widetilde{\mathcal{F}}(i, w^\alpha), [X_j^\alpha(\beta) = w_j^\alpha(\beta)]]$$

   such that

$$\widehat{T}_i^\alpha(w^\alpha)([X_j^\alpha(\beta) = w_j^\alpha(\beta)]) = \tfrac{1}{2}.$$

   Finally, by Lemma 3, extend $\widehat{T}_i^\alpha(w^\alpha)$ to a finitely additive probability measure $T_i^\alpha(w^\alpha)$ on $\mathcal{F}(i, w^\alpha)$. Define now $T_i^\alpha(u^\alpha) := T_i^\alpha(w^\alpha)$, for all $u^\alpha \in P_i(w^\alpha)$. It is easy to check that $T_i^\alpha(u^\alpha)$ and $\mathcal{F}(i, u^\alpha)$ satisfy the conditions 1–7 of the Induction hypothesis.

*Step $\alpha = \lambda$, $\lambda$ limit ordinal.* For each $P_i(u^\alpha)$ such that $u^\alpha \in W^\alpha$, choose a representing element $w^\alpha \in P_i(u^\alpha) = P_i(w^\alpha)$.

Let $T_i^{<\alpha}(w^\alpha)$ be the finitely additive probability measure defined on the field

$$\bigcup_{\beta < \alpha} \pi_{\beta,\alpha}^{-1}(\mathcal{F}(i, w^\alpha\lceil\beta)) = \bigcup_{\beta < \alpha} \pi_{\beta,\alpha}^{-1}(\mathrm{Pow}(W^\beta))$$

which is induced by $(T_i^\beta(w^\alpha\lceil\beta))_{1\le\beta<\alpha}$ (as defined in Lemma 1).

Let $\beta<\alpha$ and $E^\beta\subseteq W^\beta$ such that $\pi_{\beta,\alpha}^{-1}(E^\beta)\supseteq P_i(w^\alpha)$. By Lemma 17, we have

$$\pi_{\beta,\alpha}^{-1}(E^\beta)\supseteq\pi_{\beta+1,\alpha}^{-1}(P_i(w^\alpha\lceil\beta+1)).$$

(Note that $\beta<\alpha$ implies $\beta+1<\alpha$.) Since, by the definition of $T_i^{<\alpha}(w^\alpha)$,

$$T_i^{<\alpha}(w^\alpha)(\pi_{\beta+1,\alpha}^{-1}(P_i(w^\alpha\lceil\beta+1)))=T_i^{\beta+1}(w^\alpha\lceil\beta+1)(P_i(w^\alpha\lceil\beta+1))=1,$$

the outer measure $T_i^{<\alpha}(w^\alpha)^*(P_i(w^\alpha))$ is equal to 1. By the Łoś–Marczewski theorem, we can extend $T_i^{<\alpha}(w^\alpha)$ to a finitely additive probability measure $T_i^\alpha(w^\alpha)$ on the field

$$\mathcal{F}(i,w^\alpha)=\left[\bigcup_{\beta<\alpha}\pi_{\beta,\alpha}^{-1}(\mathrm{Pow}(W^\beta)),P_i(w^\alpha)\right]$$

such that $T_i^\alpha(w^\alpha)(P_i(w^\alpha))=1$. For $u^\alpha\in P_i(w^\alpha)$ define $T_i^\alpha(u^\alpha):=T_i^\alpha(w^\alpha)$. It is easy to check that $T_i^\alpha(u^\alpha)$ and $\mathcal{F}(i,u^\alpha)$ satisfy the conditions 1–7 of the Induction hypothesis.

*Step $\alpha=\lambda+1$, $\lambda$ limit ordinal.* For each $P_i(u^\alpha)$ such that $u^\alpha\in W^\alpha$, choose a representing element $w^\alpha\in P_i(u^\alpha)=P_i(w^\alpha)$.

Let $T_i^{<\alpha}(w^\alpha)$ be the finitely additive probability measure defined on the field

$$\pi_{\lambda,\alpha}^{-1}(\mathcal{F}(i,w^\alpha\lceil\lambda))=\left[\bigcup_{\beta<\lambda}\pi_{\beta,\alpha}^{-1}(\mathrm{Pow}(W^\beta)),\pi_{\lambda,\alpha}^{-1}(P_i(w^\alpha\lceil\lambda))\right],$$

which is induced by $T_i^\lambda(w^\alpha\lceil\lambda)$ (as defined in Lemma 1). According to Lemma 16 and the Induction hypothesis, we have for the outer measure of $P_i(w^\alpha):T_i^{<\alpha}(w^\alpha)^*(P_i(w^\alpha))=1$. So, by the Łoś–Marczewski theorem, we can extend $T_i^{<\alpha}(w^\alpha)$ to a finitely additive probability measure $\widetilde{T}_i^\alpha(w^\alpha)$ defined on the field

$$\widetilde{\mathcal{F}}(i,w^\alpha):=\left[\left[\bigcup_{\beta<\lambda}\pi_{\beta,\alpha}^{-1}(\mathrm{Pow}(W^\beta)),\pi_{\lambda,\alpha}^{-1}(P_i(w^\alpha\lceil\lambda))\right],P_i(w^\alpha)\right],$$

such that $\widetilde{T}_i^\alpha(w^\alpha)(P_i(w^\alpha))=1$.

1. Case: $w_i^\alpha(\lambda)=1$. Then

$$\pi_{\lambda,\alpha}^{-1}([\lambda-\mathrm{par}(X_j^\lambda)=\lambda-\mathrm{par}(w_j^\alpha\lceil\lambda)])=[\lambda-\mathrm{par}(X_j^\alpha)=\lambda-\mathrm{par}(w_j^\alpha)]$$
$$\supseteq P_i(w^\alpha).$$

By Lemma 3, extend $\widetilde{T}_i^\alpha(w^\alpha)$ to a finitely additive probability measure $T_i^\alpha(w^\alpha)$ on the field

$$[\pi_{\lambda,\alpha}^{-1}(\mathrm{Pow}(W^\lambda)), P_i(w^\alpha)] = \mathcal{F}(i, w^\alpha).$$

By the above, we have

$$T_i^\alpha(w^\alpha)([\lambda - \mathrm{par}(X_j^\alpha) = \lambda - \mathrm{par}(w_j^\alpha)]) = 1.$$

Define now $T_i^\alpha(u^\alpha) := T_i^\alpha(w^\alpha)$, for all $u^\alpha \in P_i(w^\alpha)$. [Note that $\mathcal{F}(i, u^\alpha) = \mathcal{F}(i, w^\alpha)$, $u_i^\alpha(\lambda) = w_i^\alpha(\lambda) = 1$ and hence, $\lambda - \mathrm{par}(w_j^\alpha) = \lambda - \mathrm{par}(u_j^\alpha)$.] It is now easy to check that $T_i^\alpha(u^\alpha)$ and $\mathcal{F}(i, u^\alpha)$ satisfy the conditions 1–7 of the Induction hypothesis.

2. Case: $w_i^\alpha(\lambda) = 0$. By Lemma 14, we have for the outer measure of $[\lambda - \mathrm{par}(X_j^\alpha) = \lambda - \mathrm{par}(w_j^\alpha)]$:

$$\widetilde{T}_i^\alpha(w^\alpha)^*([\lambda - \mathrm{par}(X_j^\alpha) = \lambda - \mathrm{par}(w_j^\alpha)]) = 1,$$

and for the inner measure

$$\widetilde{T}_i^\alpha(w^\alpha)_*([\lambda - \mathrm{par}(X_j^\alpha) = \lambda - \mathrm{par}(w_j^\alpha)]) = 0.$$

By the Łoś–Marczewski theorem, we can extend $\widetilde{T}_i^\alpha(w^\alpha)$ to a finitely additive probability measure $\widehat{T}_i^\alpha(w^\alpha)$ on the field

$$[\widetilde{\mathcal{F}}(i, w^\alpha), [\lambda - \mathrm{par}(X_j^\alpha) = \lambda - \mathrm{par}(w_j^\alpha)]]$$

such that

$$\widehat{T}_i^\alpha(w^\alpha)([\lambda - \mathrm{par}(X_j^\alpha) = \lambda - \mathrm{par}(w_j^\alpha)]) = \tfrac{1}{2}.$$

Finally, by Lemma 3, extend $\widehat{T}_i^\alpha(w^\alpha)$ to a finitely additive probability measure $T_i^\alpha(w^\alpha)$ on $\mathcal{F}(i, w^\alpha)$. Define now $T_i^\alpha(u^\alpha) := T_i^\alpha(w^\alpha)$, for all $u^\alpha \in P_i(w^\alpha)$. It is easy to check that $T_i^\alpha(u^\alpha)$ and $\mathcal{F}(i, u^\alpha)$ satisfy the conditions 1–7 of the Induction hypothesis.

*Remaining step.* To finish the proof of Theorem 2, we have to extend $T_i^\kappa(u^\kappa)$ to a finitely additive probability measure $T_i(u^\kappa)$ defined on the field $\mathrm{Pow}(W^\kappa)$ such that $T_i(u^\kappa) = T_i(w^\kappa)$, for $u^\kappa \in P_i(w^\kappa)$.

By the inductive construction, $T_i^\kappa(u^\kappa)$ is defined on

$$\left[\bigcup_{\beta < \kappa} \pi_{\beta,\kappa}^{-1}(\mathrm{Pow}(W^\beta)), P_i(w^\kappa)\right]$$

such that the conditions 1–7 of the Induction hypothesis are satisfied for $\alpha = \kappa$. For each $P_i(u^\kappa)$ such that $u^\kappa \in W^\kappa$, choose a representing element

$$w^\kappa \in P_i(u^\kappa) = P_i(w^\kappa).$$

By Lemma 3, extend $T_i^\kappa(w^\kappa)$ to a finitely additive probability measure $T_i(w^\kappa)$ on the field $\mathrm{Pow}(W^\kappa)$ and define

$$T_i(u^\kappa) := T_i(w^\kappa),$$

for $u^\kappa \in P_i(w^\kappa)$. By construction and the Induction hypothesis for $\alpha = \kappa$, the function $T_i : W^\kappa \to \Delta(W^\kappa, \mathrm{Pow}(W^\kappa))$ has all the desired properties, and hence Theorem 2 is proved. $\square$

6.1. *Proofs of Lemmas* 10 *and* 12–17.

PROOF OF LEMMA 10. By the definition of

$$\left[ \left[ \bigcup_{\beta < \gamma} \pi_{\beta,\alpha}^{-1}(\mathrm{Pow}(W^\beta)), \pi_{\gamma,\alpha}^{-1}(P_i(w^\alpha \upharpoonright \gamma)) \right], P_i(w^\alpha) \right],$$

$E$ has the form

$$E = (((\pi_{\beta,\alpha}^{-1}(A_\beta) \cap \pi_{\gamma,\alpha}^{-1}(P_i(w^\alpha \upharpoonright \gamma)))$$
$$\cup (\pi_{\eta,\alpha}^{-1}(B_\eta) \cap (W^\alpha \setminus \pi_{\gamma,\alpha}^{-1}(P_i(w^\alpha \upharpoonright \gamma))))) \cap P_i(w^\alpha))$$
$$\cup (((\pi_{\xi,\alpha}^{-1}(C_\xi) \cap \pi_{\gamma,\alpha}^{-1}(P_i(w^\alpha \upharpoonright \gamma)))$$
$$\cup (\pi_{\zeta,\alpha}^{-1}(D_\zeta) \cap (W^\alpha \setminus \pi_{\gamma,\alpha}^{-1}(P_i(w^\alpha \upharpoonright \gamma))))) \cap (W^\alpha \setminus P_i(w^\alpha))),$$

where $\beta, \eta, \xi, \zeta < \gamma$ and $A_\beta \subseteq W^\beta, B_\eta \subseteq W^\eta, C_\xi \subseteq W^\xi, D_\zeta \subseteq W^\zeta$.

The lemma follows from the following facts: If $\eta < \beta$, then $\pi_{\eta,\beta}^{-1}(B_\eta) \subseteq W^\beta$ and $\pi_{\eta,\alpha}^{-1}(B_\eta) = \pi_{\beta,\alpha}^{-1}(\pi_{\eta,\beta}^{-1}(B_\eta))$, so we can assume without loss of generality that $\beta = \eta = \xi = \zeta$. By Remark 8, we have $P_i(w^\alpha) \subseteq \pi_{\gamma,\alpha}^{-1}(P_i(w^\alpha \upharpoonright \gamma))$. $\square$

PROOF OF LEMMA 12. Let $v^\alpha \in E$. Since $\beta, o^\lambda(v_i^\alpha), o^\lambda(v_j^\alpha) < \lambda$, there is an ordinal $\xi$ such that $\max\{\beta, o^\lambda(v_i^\alpha), o^\lambda(v_j^\alpha)\} \le \xi < \lambda$ and such that the parity of $\xi + 1$ is different from $\lambda - \mathrm{par}(v_j^\alpha)$. Define now $u^\alpha \in W^\alpha$ by

$$u_0 := v_0,$$
$$u_i^\alpha := v_i^\alpha,$$
$$u_j^\alpha(\gamma) := v_j^\alpha(\gamma) \qquad \text{for all } \gamma < \alpha \text{ with } \gamma \ne \xi,$$
$$u_j^\alpha(\xi) := 1.$$

It follows that $\lambda - \mathrm{par}(u_j^\alpha) \ne \lambda - \mathrm{par}(v_j^\alpha)$ and $u^\alpha \upharpoonright \beta = v^\alpha \upharpoonright \beta$, which implies $u^\alpha \in E$.

(a) If $v^\alpha \in P_i(w^\alpha)$, then it is easy to check that $u^\alpha \in P_i(v^\alpha) = P_i(w^\alpha)$.

(b) If $v^\alpha \in \pi_{\lambda,\alpha}^{-1}(P_i(w^\alpha \lceil \lambda)) \cap (W^\alpha \setminus P_i(w^\alpha))$, then it follows that $v^\alpha \in \pi_{\lambda,\alpha}^{-1}(P_i(w^\alpha \lceil \lambda))$ and $v_i^\alpha(\lambda) = 1$. It is again easy to check that $u^\alpha \in \pi_{\lambda,\alpha}^{-1}(P_i(v^\alpha \lceil \lambda)) = \pi_{\lambda,\alpha}^{-1}(P_i(w^\alpha \lceil \lambda))$ and since $u_i^\alpha(\lambda) = 1$, we have $u^\alpha \in (W^\alpha \setminus P_i(w^\alpha))$.

(c) If $v^\alpha \notin \pi_{\lambda,\alpha}^{-1}(P_i(w^\alpha \lceil \lambda))$, then there are four cases:

1. $v_i^\alpha \lceil \lambda \neq w_i^\alpha \lceil \lambda$. From $u_i^\alpha \lceil \lambda = v_i^\alpha \lceil \lambda$ it follows that $u^\alpha \notin \pi_{\lambda,\alpha}^{-1}(P_i(w^\alpha \lceil \lambda))$.

2. There is a $\gamma < \lambda$ such that
$$(v_i^\alpha \lceil \lambda)(\gamma + 1) = (w_i^\alpha \lceil \lambda)(\gamma + 1) = 1 \quad \text{and} \quad (v_j^\alpha \lceil \lambda)(\gamma) \neq (w_j^\alpha \lceil \lambda)(\gamma).$$

   Since $\gamma + 1 < \lambda$ and $\max\{\beta, o^\lambda(v_i^\alpha), o^\lambda(v_j^\alpha)\} > \gamma + 1$, it follows that $u^\alpha \lceil \gamma + 2 = v^\alpha \lceil \gamma + 2$ and therefore $u^\alpha \notin \pi_{\lambda,\alpha}^{-1}(P_i(w^\alpha \lceil \lambda))$.

3. There is a limit ordinal $\widehat{\lambda} < \lambda$ such that
$$(v_i^\alpha \lceil \lambda)(\widehat{\lambda}) = (w_i^\alpha \lceil \lambda)(\widehat{\lambda}) = 1 \quad \text{and} \quad \widehat{\lambda} - \mathrm{par}(v_j^\alpha) \neq \widehat{\lambda} - \mathrm{par}(w_j^\alpha).$$

   We have $\xi > \widehat{\lambda}$, and therefore
$$(u_i^\alpha \lceil \lambda)(\widehat{\lambda}) = (v_i^\alpha \lceil \lambda)(\widehat{\lambda}) \quad \text{and} \quad \widehat{\lambda} - \mathrm{par}(u_j^\alpha) = \widehat{\lambda} - \mathrm{par}(v_j^\alpha).$$

   It follows that $u^\alpha \notin \pi_{\lambda,\alpha}^{-1}(P_i(w^\alpha \lceil \lambda))$.

4. $(v_i^\alpha \lceil \lambda)(0) = (w_i^\alpha \lceil \lambda)(0) = 1$ and $v_0 \neq w_0$. We have
$$(u_i^\alpha \lceil \lambda)(0) = (v_i^\alpha \lceil \lambda)(0) = 1 \quad \text{and} \quad u_0 = v_0.$$

   It follows that $u^\alpha \notin \pi_{\lambda,\alpha}^{-1}(P_i(w^\alpha \lceil \lambda))$. □

PROOF OF LEMMA 13. Let $v^\alpha \in E$. Define $u^\alpha \in W^\alpha$ by
$$u_0 := v_0,$$
$$u_i^\alpha := v_i^\alpha,$$
$$u_j^\alpha(\gamma) := v_j^\alpha(\gamma) \qquad \text{for all } \gamma < \alpha \text{ with } \gamma \neq \beta,$$
$$u_j^\alpha(\beta) := 1 - v_j^\alpha(\beta).$$

It follows that $u_j^\alpha(\beta) \neq v_j^\alpha(\beta)$ and $u^\alpha \lceil \beta = v^\alpha \lceil \beta$, which implies $u^\alpha \in E$.

(a) If $v^\alpha \in P_i(w^\alpha)$, then it is easy to check that $u^\alpha \in P_i(v^\alpha) = P_i(w^\alpha)$.

(b) If $v^\alpha \in \pi_{\beta+1,\alpha}^{-1}(P_i(w^\alpha \lceil \beta + 1)) \cap (W^\alpha \setminus P_i(w^\alpha))$, then it follows that $v^\alpha \in \pi_{\beta+1,\alpha}^{-1}(P_i(w^\alpha \lceil \beta + 1))$ and $v_i^\alpha(\beta + 1) = 1$. It is again easy to check that $u^\alpha \in \pi_{\beta+1,\alpha}^{-1}(P_i(v^\alpha \lceil \beta + 1)) = \pi_{\beta+1,\alpha}^{-1}(P_i(w^\alpha \lceil \beta + 1))$ and since $u_i^\alpha(\beta + 1) = 1$, we have $u^\alpha \in (W^\alpha \setminus P_i(w^\alpha))$.

(c) If $v^\alpha \notin \pi_{\beta+1,\alpha}^{-1}(P_i(w^\alpha \lceil \beta + 1))$, then there are four cases:

1. $v_i^\alpha \lceil \beta + 1 \neq w_i^\alpha \lceil \beta + 1$. From $u_i^\alpha \lceil \beta + 1 = v_i^\alpha \lceil \beta + 1$, it follows that $u^\alpha \notin \pi_{\beta+1,\alpha}^{-1}(P_i(w^\alpha \lceil \beta + 1))$.

2. There is a $\gamma < \beta$ such that

$$(v_i^\alpha \lceil \beta + 1)(\gamma + 1) = (w_i^\alpha \lceil \beta + 1)(\gamma + 1) = 1$$

and

$$(v_j^\alpha \lceil \beta + 1)(\gamma) \neq (w_j^\alpha \lceil \beta + 1)(\gamma).$$

By the definition of $u^\alpha$, $(u_i^\alpha \lceil \beta + 1)(\gamma + 1) = 1$ and, since $\gamma < \beta$,

$$(u_j^\alpha \lceil \beta + 1)(\gamma) = (u_j^\alpha \lceil \beta + 1)(\gamma),$$

hence $u^\alpha \notin \pi_{\beta+1,\alpha}^{-1}(P_i(w^\alpha \lceil \beta + 1))$.

3. There is a limit ordinal $\lambda < \beta + 1$ such that

$$(v_i^\alpha \lceil \beta + 1)(\lambda) = (w_i^\alpha \lceil \beta + 1)(\lambda) = 1 \quad \text{and} \quad \lambda - \mathrm{par}(v_j^\alpha) \neq \lambda - \mathrm{par}(w_j^\alpha).$$

We have $(u_i^\alpha \lceil \beta + 1)(\lambda) = 1$ and, since $\lambda \leq \beta$, $\lambda - \mathrm{par}(u_j^\alpha) = \lambda - \mathrm{par}(v_j^\alpha)$. It follows that $u^\alpha \notin \pi_{\beta+1,\alpha}^{-1}(P_i(w^\alpha \lceil \beta + 1))$.

4. $(v_i^\alpha \lceil \beta + 1)(0) = (w_i^\alpha \lceil \beta + 1)(0) = 1$ and $v_0 \neq w_0$. We have

$$(u_i^\alpha \lceil \beta + 1)(0) = (v_i^\alpha \lceil \beta + 1)(0) = 1 \quad \text{and} \quad u_0 = v_0.$$

It follows that $u^\alpha \notin \pi_{\beta+1,\alpha}^{-1}(P_i(w^\alpha \lceil \beta + 1))$. $\quad\square$

PROOF OF LEMMA 14. The first point of the lemma follows from Lemmas 10 and 12.

The second point follows from the first and the fact that

$$\left[ \left[ \bigcup_{\beta < \lambda} \pi_{\beta,\alpha}^{-1}(\mathrm{Pow}(W^\beta)), \pi_{\lambda,\alpha}^{-1}(P_i(w^\alpha \lceil \lambda)) \right], P_i(w^\alpha) \right]$$

is a field on $W^\alpha$ (and therefore it is closed under complements).

The last two points of the lemma follow directly from the first two points. $\square$

PROOF OF LEMMA 15. Note that if $\alpha \geq \beta + 1$, then

$$\bigcup_{\xi < \beta + 1} \pi_{\xi,\alpha}^{-1}(\mathrm{Pow}(W^\xi)) = \pi_{\beta,\alpha}^{-1}(\mathrm{Pow}(W^\beta)).$$

The proof is now analogous to the proof of Lemma 14—just replace $\lambda$ by $\beta + 1$ and Lemma 12 by Lemma 13. $\quad\square$

PROOF OF LEMMA 16. Since $P_i(w^\alpha) \subseteq \pi_{\gamma,\alpha}^{-1}(P_i(w^\alpha \lceil \gamma))$, it follows from the definition of

$$\left[ \bigcup_{\beta < \gamma} \pi_{\beta,\alpha}^{-1}(\mathrm{Pow}(W^\beta)), \pi_{\gamma,\alpha}^{-1}(P_i(w^\alpha \lceil \gamma)) \right],$$

that there is a $\beta < \gamma$ and an $E^\beta \subseteq W^\beta$ such that

$$E \cap \pi_{\gamma,\alpha}^{-1}(P_i(w^\alpha \lceil \gamma)) = \pi_{\beta,\alpha}^{-1}(E^\beta) \cap \pi_{\gamma,\alpha}^{-1}(P_i(w^\alpha \lceil \gamma)) \supseteq P_i(w^\alpha).$$

*Claim.* $\pi_{\beta,\alpha}^{-1}(E^\beta) \supseteq \pi_{\gamma,\alpha}^{-1}(P_i(w^\alpha \lceil \gamma)).$

Assume to the contrary, that there is a $v^\alpha \in \pi_{\gamma,\alpha}^{-1}(P_i(w^\alpha \lceil \gamma)) \setminus \pi_{\beta,\alpha}^{-1}(E^\beta)$. Since $\beta + 1 \le \gamma$, we have $v^\alpha \lceil \beta + 1 \in P_i(w^\alpha \lceil \beta + 1)$. By Lemma 11, there is a $u^\alpha \in P_i(w^\alpha)$ such that $u^\alpha \lceil \beta = v^\alpha \lceil \beta$. Since $E^\beta \subseteq W^\beta$ and $v^\alpha \notin \pi_{\beta,\alpha}^{-1}(E^\beta)$, it follows that $u^\alpha \notin \pi_{\beta,\alpha}^{-1}(E^\beta)$, a contradiction to $\pi_{\beta,\alpha}^{-1}(E^\beta) \supseteq P_i(w^\alpha)$.  $\square$

PROOF OF LEMMA 17.   Assume that there is a $v^\lambda \in \pi_{\beta+1,\lambda}^{-1}(P_i(w^\lambda \lceil \beta + 1)) \setminus \pi_{\beta,\lambda}^{-1}(E^\beta)$. By Lemma 11, there is a $u^\lambda \in P_i(w^\lambda)$ such that $u^\lambda \lceil \beta = v^\lambda \lceil \beta$. Therefore $u^\lambda \notin \pi_{\beta,\lambda}^{-1}(E^\beta)$, a contradiction to $\pi_{\beta,\lambda}^{-1}(E^\beta) \supseteq P_i(w^\lambda)$.  $\square$

## APPENDIX

PROOF OF LEMMA 7.   We prove the lemma by induction on the formation of $\kappa$-expressions.

(a) Let $\varphi = E$, where $E \in \Sigma_{\{h,t\}} = \text{Pow}(\{h,t\})$, and let $w^\kappa, u^\kappa \in W^\kappa$ such that $w^\kappa \lceil 0 = u^\kappa \lceil 0$. By definition, $v^\kappa \in E^{W^\kappa}$ iff $\theta(v^\kappa) \in E$, for $v^\kappa \in W^\kappa$. But we have $\theta(v^\kappa) = v_0 = v^\kappa \lceil 0$, for $v^\kappa \in W^\kappa$. It follows that $u^\kappa \in E^{W^\kappa}$ iff $w^\kappa \in E^{W^\kappa}$.

(b) Let $\varphi = \neg \psi$ such that $\text{depth}(\varphi) \le \alpha$ and let $w^\kappa, u^\kappa \in W^\kappa$ such that $w^\kappa \lceil \alpha = u^\kappa \lceil \alpha$. It follows that $\text{depth}(\psi) \le \alpha$ and $u^\kappa \in \varphi^{W^\kappa}$ iff $u^\kappa \notin \psi^{W^\kappa}$ iff—by the induction assumption—$w^\kappa \notin \psi^{W^\kappa}$, which is the case iff $w^\kappa \in \varphi^{W^\kappa}$.

(c) Let $p \in [0,1]$, $i \in \{a,b\}$, $\varphi = B_i^p(\psi)$, $\text{dp}(\psi) + 1 = \beta + 1 \le \alpha$ and $w^\kappa, u^\kappa \in W^\kappa$ such that $w^\kappa \lceil \alpha = u^\kappa \lceil \alpha$. By the induction assumption, there is a $E^\beta \subseteq W^\beta$ such that $\psi^{W^\kappa} = \pi_{\beta,\kappa}^{-1}(E^\beta)$. By Theorem 2 and Remark 11,

$$T_i(u^\kappa)(\psi^{W^\kappa}) = T_i(w^\kappa)(\psi^{W^\kappa}).$$

It follows that $u^\kappa \in (B_i^p(\psi))^{W^\kappa}$ iff $w^\kappa \in (B_i^p(\psi))^{W^\kappa}$.

(d) Let $|\Psi| < \kappa$, $\varphi = \bigwedge_{\psi \in \Psi} \psi$ such that $\text{depth}(\varphi) \le \alpha$, and let $w^\kappa, u^\kappa \in W^\kappa$ such that $w^\kappa \lceil \alpha = u^\kappa \lceil \alpha$. Then $\text{depth}(\psi) \le \alpha$, for $\psi \in \Psi$. By the induction assumption, $u^\kappa \in \psi^{W^\kappa}$ iff $w^\kappa \in \psi^{W^\kappa}$, for $\psi \in \Psi$. It follows that $u^\kappa \in \varphi^{W^\kappa}$ iff $w^\kappa \in \varphi^{W^\kappa}$.
$\square$

Lemma 8 follows directly from Theorem 2 and Remark 11.

PROOF OF LEMMA 9.   The first point is clear. We show the second and the third points by a transfinite induction on $0 \le \beta < \kappa$.

(a) $\beta = 0$. According to Lemma 8 and the first point of this lemma, we have

$$(\varphi_i^1(0))^{W^\kappa} := (B_i^1(\{h\}) \vee B_i^1(\{t\}))^{W^\kappa} = [X_i^\kappa(0) = 1]$$

and

$$(\varphi_i^0(0))^{W^\kappa} := (\neg(B_i^1(\{h\}) \vee B_i^1(\{t\})))^{W^\kappa} = [X_i^\kappa(0) = 0].$$

And obviously, $\mathrm{dp}(\varphi_i^0(0)) = \mathrm{dp}(\varphi_i^1(0)) = 1$.

(b) $\beta = \gamma + 1$. According to the induction assumption, there are $\kappa$-expressions $\varphi_j^0(\gamma)$ and $\varphi_j^1(\gamma)$ such that

$$(\varphi_j^0(\gamma))^{W^\kappa} = [X_j^\kappa(\gamma) = 0],$$

$$(\varphi_j^1(\gamma))^{W^\kappa} = [X_j^\kappa(\gamma) = 1],$$

$$\mathrm{dp}(\varphi_j^0(\gamma)) = \mathrm{dp}(\varphi_j^1(\gamma))$$

$$= \beta.$$

Define

$$\varphi_i^1(\beta) := B_i^1(\varphi_j^0(\gamma)) \vee B_i^1(\varphi_j^1(\gamma))$$

and

$$\varphi_i^0(\beta) := \neg\varphi_i^1(\beta).$$

$\varphi_i^1(\beta)$ and $\varphi_i^0(\beta)$ are $\kappa$-expressions of depth $\beta + 1$. We have

$$[X_i^\kappa(\beta) = 0] = W^\kappa \setminus [X_i^\kappa(\beta) = 1].$$

By Lemma 8 and the induction assumption, it follows that

$$[X_i^\kappa(\beta) = 1] = (\varphi_i^1(\beta))^{W^\kappa}$$

and

$$[X_i^\kappa(\beta) = 0] = (\varphi_i^0(\beta))^{W^\kappa}.$$

(c) Let $\lambda < \kappa$ be a limit ordinal. For $i \in \{a, b\}$ and $\beta < \lambda$ define in $\underline{W}^\kappa$:

$$[Y_i^\lambda(\beta)] := [X_i^\kappa(\beta) = 1] \cap \bigcap_{\beta < \alpha < \lambda} [X_i^\kappa(\alpha) = 0],$$

$$[Z_i^\lambda] := \bigcap_{0 \le \alpha < \lambda} [X_i^\kappa(\alpha) = 0].$$

According to the induction assumption for $\alpha, \beta < \lambda$ and the fact that $|\lambda| < \kappa$, it follows that

$$\psi_i^\lambda(\beta) := \varphi_i^1(\beta) \wedge \bigwedge_{\beta < \alpha < \lambda} \varphi_i^0(\alpha)$$

and

$$\chi_i^\lambda := \bigwedge_{0 \le \alpha < \lambda} \varphi_i^0(\alpha)$$

are $\kappa$-expressions such that

$$\mathrm{dp}(\psi_i^\lambda(\beta)) = \max\{\beta + 1, \sup\{\mathrm{dp}(\varphi_i^0(\alpha))|\beta < \alpha < \lambda\}\}$$
$$= \sup\{\alpha + 1|\alpha < \lambda\}$$
$$= \lambda,$$

and, similarly, $\mathrm{dp}(\chi_i^\lambda) = \lambda$.

It follows from the induction assumption that

$$[Y_i^\lambda(\beta)] = (\psi_i^\lambda(\beta))^{W^\kappa} \quad \text{and} \quad [Z_i^\lambda] = (\chi_i^\lambda)^{W^\kappa}.$$

Since $o^\lambda(w_i^\kappa)$, for $w_i^\kappa \in W_i^\kappa$, can never be a limit ordinal, we have

$$[\lambda - \mathrm{par}(X_i^\kappa) = \mathrm{even}] = [Z_i^\lambda] \cup \bigcup_{\beta < \lambda, \beta \, \mathrm{odd}} [Y_i^\lambda(\beta)]$$

and

$$[\lambda - \mathrm{par}(X_i^\kappa) = \mathrm{odd}] = \bigcup_{\beta < \lambda, \beta \, \mathrm{even}} [Y_i^\lambda(\beta)].$$

Again, since $|\lambda| < \kappa$, it follows from the above that

$$\varphi_i^{\mathrm{even}}(\lambda) := \chi_i^\lambda \vee \bigvee_{\beta < \lambda, \beta \, \mathrm{odd}} \psi_i^\lambda(\beta)$$

and

$$\varphi_i^{\mathrm{odd}}(\lambda) := \bigvee_{\beta < \lambda, \beta \, \mathrm{even}} \psi_i^\lambda(\beta)$$

are $\kappa$-expressions such that

$$\mathrm{dp}(\varphi_i^{\mathrm{even}}(\lambda)) = \max\{\mathrm{dp}(\chi_i^\lambda), \sup\{\mathrm{dp}(\psi_i^\lambda(\beta))|\beta < \lambda, \beta \, \mathrm{odd}\}\} = \lambda,$$

and

$$\mathrm{dp}(\varphi_i^{\mathrm{odd}}(\lambda)) = \lambda.$$

By the definitions and the induction assumption, we have

$$(\varphi_i^{\mathrm{even}}(\lambda))^{W^\kappa} = [\lambda - \mathrm{par}(X_i^\kappa) = \mathrm{even}]$$

and

$$(\varphi_i^{\mathrm{odd}}(\lambda))^{W^\kappa} = [\lambda - \mathrm{par}(X_i^\kappa) = \mathrm{odd}].$$

(d) $\beta = \lambda$, $\lambda$ limit ordinal $< \kappa$. By Lemma 8 and the above we have

$$(\varphi_i^1(\lambda))^{W^\kappa} := (B_i^1(\varphi_j^{\text{even}}(\lambda)) \vee B_i^1(\varphi_j^{\text{odd}}(\lambda)))^{W^\kappa} = [X_i^\kappa(\lambda) = 1],$$

$$(\varphi_i^0(\lambda))^{W^\kappa} := (\neg(B_i^1(\varphi_j^{\text{even}}(\lambda)) \vee B_i^1(\varphi_j^{\text{odd}}(\lambda))))^{W^\kappa} = [X_i^\kappa(\lambda) = 0]$$

and

$$\begin{aligned}
&\operatorname{dp}(B_i^1(\varphi_j^{\text{even}}(\lambda)) \vee B_i^1(\varphi_j^{\text{odd}}(\lambda))) \\
&= \operatorname{dp}(\neg(B_i^1(\varphi_j^{\text{even}}(\lambda)) \vee B_i^1(\varphi_j^{\text{odd}}(\lambda))))^{W^\kappa} \\
&= \lambda + 1. \qquad \qquad \qquad \square
\end{aligned}$$

**Acknowledgments.** I would like to thank an anonymous referee for very valuable suggestions that improved the readability of the paper a lot. Helpful comments of Aviad Heifetz, Jean-François Mertens and Dov Samet are gratefully acknowledged.

Instituto de Análisis Económico, CSIC
Campus UAB
08193 Bellaterra, Barcelona
Spain
e-mail: Martin.Meier@uab.es